\documentclass[11pt]{article}
\usepackage{amsmath}
\usepackage{a4wide,amsfonts,amsmath,latexsym,amssymb,euscript,graphicx,units,mathrsfs}
\usepackage{graphicx}
\usepackage[T1]{fontenc}
\usepackage{latexsym}
\usepackage{xypic}
\usepackage{eufrak}
\usepackage{euscript}
\usepackage{verbatim}
\usepackage{fancyhdr}
\usepackage[english]{babel}
\usepackage{mathrsfs}
\usepackage{units}

\newtheorem{theorem}{Theorem}

\newtheorem{corollary}[theorem]{Corollary}

\newtheorem{definition}[theorem]{Definition}

\newtheorem{lemma}[theorem]{Lemma}

\newtheorem{proposition}[theorem]{Proposition}

\newenvironment{proof}[1][Proof]{\textbf{#1.} }{\ \rule{0.5em}{0.5em}}

\newcommand{\R}{\mathbb{R}}

\def\qed{ \hfill \vrule width.25cm height.25cm depth0cm\smallskip}
\def\HH{{\mathcal H}}
\begin{document}

\renewcommand{\thefootnote}{\fnsymbol{footnote}}

\begin{center}
{\large \textbf{Central and non-central limit theorems for weighted
power
variations of fractional Brownian motion}}\\[0pt]
~\\[0pt]
Ivan Nourdin\footnote{%
Laboratoire de Probabilit{\'e}s et Mod{\`e}les Al{\'e}atoires, Universit{\'e}
Pierre et Marie Curie, Bo{\^\i}te courrier 188, 4 Place Jussieu, 75252 Paris
Cedex 5, France, \texttt{ivan.nourdin@upmc.fr}}, David Nualart\footnote{%
Department of Mathematics, University of Kansas, 405 Snow Hall, Lawrence,
Kansas 66045-2142, USA, \texttt{nualart@math.ku.edu}}\footnote{The work of D. Nualart is supported by the NSF Grant DMS-0604207} and Ciprian A. Tudor%
\footnote{%
SAMOS/MATISSE, Centre d'\'Economie de La Sorbonne, Universit\'e de
Panth\'eon-Sorbonne Paris 1, 90 rue de Tolbiac, 75634 Paris Cedex 13,
France, \texttt{tudor@univ-paris1.fr}}\\[0pt]
~\\[0pt]
\end{center}

{\small \noindent \textbf{Abstract: } In this paper, we prove some
central and non-central limit theorems for renormalized weighted
power variations of order $q\ge 2 $ of the fractional Brownian
motion with Hurst parameter $H\in (0,1)$, where $q$ is an integer.
The central limit holds for $\frac 1{2q} <H\leq 1-\frac 1{2q} $,
the limit being a conditionally Gaussian distribution. If $%
H<\frac 1{2q}$ we show the convergence in $L^2$ to a limit which only
depends on the fractional Brownian motion, and if $H> 1-\frac 1{2q}$ we show
the convergence in $L^2$ to a stochastic integral with respect to the Hermite
process of order $q $. \newline
}

{\small \noindent \textbf{Key words:} fractional Brownian motion, central
limit theorem, non-central limit theorem, Hermite process.~\newline
}

{\small \noindent{\bf 2000 Mathematics Subject Classification:}
60F05, 60H05, 60G15, 60H07. ~\newline}

{\small \noindent{\it This version}: August 2009.}

\section{Introduction}

 \setcounter{footnote}{0}

The study of single path behavior of stochastic processes is often based on
the study of their power variations, and there exists a very extensive
literature on the subject. Recall that, a real $q>0$ being given, the $%
q $-power variation of a stochastic process $X$, with respect to a
subdivision $\pi _{n}=\{0=t_{n,0}<t_{n,1}<\ldots <t_{n,
\kappa(n)}=1\}$ of $[0,1]$,
is defined to be the sum
\begin{equation*}
%%% debut
\sum_{k=1}^{\kappa(n)}|X_{t_{n,k}}-X_{t_{n,k-1}}|^{q}.
%%% fin
\end{equation*}%
For simplicity, consider from now on the case where $t_{n,k}=k2^{-n}$
for $n\in\{ 1,2,3,\ldots\}$ and $k\in \{0,\ldots ,2^{n}\}$. In the present paper we
wish to point out some interesting phenomena when $X=B$ is a fractional
Brownian motion of Hurst index $H\in (0,1)$, and when $q\geq 2$ is an
integer. In fact, we will also drop the absolute value (when $q$ is odd) and
we will introduce some weights. More precisely, we will consider
\begin{equation}
\sum_{k=1}^{2^{n}}f(B_{(k-1)2^{-n}})(\Delta B_{k2^{-n}})^{q},\quad q\in
\{2,3,4,\ldots \},  \label{sum-intro}
\end{equation}%
where the function $f:\mathbb{R}\rightarrow \mathbb{R}$ is assumed to be
smooth enough and where $\Delta B_{k2^{-n}}$ denotes, here and in all the paper, the increment $%
B_{k2^{-n}}-B_{(k-1)2^{-n}}$.

The analysis of the asymptotic behavior of quantities of type (\ref%
{sum-intro}) is motivated, for instance, by the study of the exact
rates of convergence of some approximation schemes of scalar
stochastic differential equations driven by $B$ (see \cite{GN},
\cite{NN} and \cite{N06}) besides, of course, the traditional
applications of quadratic variations to parameter estimation
problems.

Now, let us recall some known results concerning $q$-power variations (for $%
q=2,3,4,\ldots $), which are today more or less classical. First, assume
that the Hurst index is $H=\frac{1}{2}$, that is $B$ is a standard Brownian
motion. Let $\mu _{q}$ denote the $q$th moment of a standard Gaussian random
variable $G\sim \mathscr{N}(0,1)$. By the scaling property of the Brownian
motion and using the central limit theorem, it is immediate that, as $%
n\rightarrow \infty $:
\begin{equation}
2^{-n/2}\sum_{k=1}^{2^{n}}\left[ (2^{n/2}\Delta B_{k2^{-n}})^{q}-\mu _{q}%
\right] \,\,\,{\overset{\mathrm{Law}}{\longrightarrow }}\,\,\,\mathscr{N}%
(0,\mu _{2q}-\mu _{q}^{2}).  \label{mbpa}
\end{equation}%
When weights are introduced, an interesting phenomenon appears: instead of
Gaussian random variables, we rather obtain mixing random variables as limit
in (\ref{mbpa}). Indeed, when $q$ is even
%%% debut
and $f:\R\to\R$ is continuous and has polynomial growth,
%%% fin
it is a very particular case of a
more general result by Jacod \cite{Jacod}
%%% debut
(see also Section 2 in Nourdin and Peccati \cite{NoPe} for related results)
%%% fin
that we have, as $n\rightarrow
\infty $:
\begin{equation}
2^{-n/2}\sum_{k=1}^{2^{n}}f(B_{(k-1)2^{-n}})\left[ (2^{n/2}\Delta
B_{k2^{-n}})^{q}-\mu _{q}\right] \,\,\,{\overset{\mathrm{Law}}{%
\longrightarrow }}\,\,\,\sqrt{\mu _{2q }-\mu _{q}^{2}}%
\int_{0}^{1}f(B_{s})dW_{s}.  \label{pombpa}
\end{equation}%
Here, $W$ denotes another standard Brownian motion, independent of $B$. When
$q$ is odd,
%%% debut
still for $f:\R\to\R$ continuous with polynomial growth,
%%% fin
we have, this time, as $n\rightarrow \infty $:
\begin{equation}
2^{-n/2}\sum_{k=1}^{2^{n}}f(B_{(k-1)2^{-n}})(2^{n/2}\Delta
B_{k2^{-n}})^{q}\,\,\,{\overset{\mathrm{Law}}{\longrightarrow }}%
\,\,\,\int_{0}^{1}f(B_{s})\big(\sqrt{\mu _{2q}-\mu _{q+1}^{2}}\,dW_{s}+\mu
_{q+1}\,dB_{s}\big),  \label{pombimpa}
\end{equation}%
see for instance \cite{NoPe}.

Secondly, assume that $H\neq\frac{1}{2}$, that is the case where the
fractional Brownian motion $B$ has not independent increments anymore. Then (%
\ref{mbpa}) has been extended by Breuer and Major \cite{BM}, Dobrushin and
Major \cite{DM}, Giraitis and Surgailis \cite{GS} or Taqqu \cite{Taqqu}.
Precisely, five cases are considered, according to the evenness of $q$ and
the value of $H$:

\begin{itemize}
\item[\textbullet] if $q$ is even and if $H\in(0,\frac{3}{4})$, as $%
n\to\infty$,
\begin{equation}  \label{mbfpa}
2^{-n/2} \sum_{k=1}^{2^n}\big[(2^{nH}\Delta B_{k2^{-n}})^q-\mu_{q}\big] %
\,\,\,{\overset{\mathrm{Law}}{\longrightarrow}}\,\,\, \mathscr{N}%
(0,\widetilde{\sigma}^2_{H,q}).
\end{equation}

\item[\textbullet] if $q$ is even and if $H=\frac{3}{4}$, as $n\to\infty$,
\begin{equation}  \label{mbfpa-2}
 \frac 1 {\sqrt{n}} 2^{-n/2} \sum_{k=1}^{2^n}\big[(2^{\frac 34 n}\Delta B_{k2^{-n}})^q-\mu_{q}\big] %
\,\,\,{\overset{\mathrm{Law}}{\longrightarrow}}\,\,\, \mathscr{N}%
(0,\widetilde{\sigma}^2_{\frac 34,q}).
\end{equation}

\item[\textbullet] if $q$ is even and if $H\in(\frac{3}{4},1)$, as $%
n\to\infty$,
\begin{equation}  \label{mbfpa2}
2^{n-2nH} \sum_{k=1}^{2^n}\big[(2^{nH}\Delta B_{k2^{-n}})^q-\mu_{q}\big] %
\,\,\,{\overset{\mathrm{Law}}{\longrightarrow}}\,\,\,
\mbox{``Hermite r.v.''.}
\end{equation}

\item[\textbullet] if $q$ is odd and if $H\in(0,\frac{1}{2}]$, as $%
n\to\infty $,
\begin{equation}  \label{mbfim}
2^{-n/2} \sum_{k=1}^{2^n}(2^{nH}\Delta B_{k2^{-n}})^q \,\,\,{\overset{%
\mathrm{Law}}{\longrightarrow}}\,\,\, \mathscr{N}(0,\widetilde{\sigma}^2_{H,q}).
\end{equation}

\item[\textbullet] if $q$ is odd and if $H\in(\frac{1}{2},1)$, as $%
n\to\infty $,
\begin{equation}  \label{mbfim2}
2^{-nH} \sum_{k=1}^{2^n}(2^{nH}\Delta B_{k2^{-n}})^q \,\,\,{\overset{\mathrm{%
Law}}{\longrightarrow}}\,\,\, \mathscr{N}(0,\widetilde{\sigma}^2_{H,q}).
\end{equation}
\end{itemize}

%%% debut
Here, $\widetilde{\sigma}_{H,q}>0$ denote some constant depending only
on $H$ and $q$.
%%% fin
The term
``Hermite r.v.'' denotes a
random variable whose distribution is the same as that of $Z^{(2)}$ at time
one, for $Z^{(2)}$ defined in Definition \ref{d1} below.

Now, let us proceed with the results concerning the \textit{weighted} power
variations in the case where $H\neq \frac{1}{2}$.
Consider the following condition on a function  $f:\mathbb{R}\rightarrow \mathbb{R}$, where $q\ge 2$ is an integer:  \begin{equation*}
(\mathbf{H}_{q})\quad
\mbox{\textit{$f$ belongs to $\mathscr{C}^{2q}$ and, for any $p\in(0,\infty)$ and $0\le i\le 2q$: $\sup_{t\in[0,1]}
E\left\{|f^{(i)}(B_t)|^p \right\}<\infty$.}}
\end{equation*}
 Suppose that $f$ satisfies  $(\mathbf{H}_{q})$.
If $q$ is even and $H\in(\frac{1}{2},\frac{3%
}{4})$, then by Theorem 2 in Le\'{o}n and Lude\~{n}a \cite{LL} (see also
Corcuera \textit{et al} \cite{CNW} for related results on the asymptotic
behavior of the $p$-variation of stochastic integrals with respect to $B$)
we have, as $n\rightarrow\infty$:
\begin{equation}  \label{pombfpa}
2^{-n/2} \sum_{k=1}^{2^n}f(B_{(k-1)2^{-n}}) \big[(2^{nH}\Delta
B_{k2^{-n}})^q-\mu_{q}\big] \,\,\,{\overset{\mathrm{Law}}{\longrightarrow}}%
\,\,\, \widetilde{\sigma}_{H,q}\int_0^1 f(B_s)dW_s,
\end{equation}
where, once again, $W$ denotes a standard Brownian motion independent of $B$
%%% debut
while $\widetilde{\sigma}_{H,q}$ is the constant appearing in (\ref{mbfpa}).
%%% fin
Thus, (\ref{pombfpa}) shows for (\ref{sum-intro}) a similar behavior to
that observed in the standard Brownian case, compare with (\ref{pombpa}). In
contradistinction, the asymptotic behavior of (\ref{sum-intro}) can be
completely different of (\ref{pombpa}) or (\ref{pombfpa}) for other values
of $H$. The first result in this direction has been observed by Gradinaru
\textit{et al} \cite{GRV}.
Namely, if $q\geq 3$ is odd
and $H\in(0,\frac12)$, we have, as $n\rightarrow\infty$:
\begin{equation}  \label{pomfbim}
2^{nH-n}\sum_{k=1}^{2^n}f(B_{(k-1)2^{-n}}) (2^{nH}\Delta B_{k2^{-n}})^q
\,\,\,{\overset{L^2}{\longrightarrow}}\,\,\, -\frac{\mu_{q+1}}{2}%
\int_0^1 f^{\prime}(B_s)ds.
\end{equation}
Also, when $q=2$ and $H\in (0,\frac{1}{4})$, Nourdin \cite{Nourdin} proved
that we have, as $n\rightarrow\infty$:
\begin{equation}
2^{2Hn-n}\sum_{k=1}^{2^n}f(B_{(k-1)2^{-n}})\big[ (2^{nH}\Delta
B_{k2^{-n}})^{2}-1\big] \,\,{\overset{L^{2}}{\longrightarrow }}\,\,%
\frac{1}{4}\int_{0}^{1}f^{\prime \prime }(B_{s})ds.  \label{cv-L2}
\end{equation}

In view of (\ref{pombpa}), (\ref{pombimpa}), (\ref{pombfpa}), (\ref{pomfbim}%
) and (\ref{cv-L2}), we observe that the asymptotic behaviors of the
power variations of fractional Brownian motion (\ref{sum-intro}) can
be really different, depending on the values of $q$ and $H$. The aim of the present paper is to
investigate what happens in the whole generality with respect to $q$
and $H$. Our main tool is the Malliavin calculus that appeared, in
several recent papers, to be very useful in the study of the power
variations for stochastic processes. As we will see, the Hermite
polynomials play a crucial role in this analysis. In the sequel, for
an integer $q\ge 2$, we write $H_q$ for the Hermite polynomial with
degree $q$ defined by
\begin{equation*}
H_q(x)=\frac{(-1)^q}{q!}e^{\frac{x^2}{2}}\frac{d^q}{dx^q}\left(%
e^{-\frac{x^2}{2}}\right),
\end{equation*}
and we consider, when $f:\mathbb{R}\rightarrow\mathbb{R}$ is a deterministic
function, the sequence of \textit{weighted Hermite variation of order $q$}
defined by
\begin{equation}  \label{weighted-hermite}
V^{(q)} _{n}(f):= \sum_{k=1}^{2^n} f\left( B _{(k-1)2^{-n}}\right) \,H_{q}%
\big( 2^{nH} \Delta B_{k2^{-n}}\big).
\end{equation}
The following is the main result of this paper.

\begin{theorem}
\label{main} Fix an integer $q\geq 2$, and suppose   that $f$ satisfies  $(\mathbf{H}_{q})$.

\begin{description}
\item[1.]  Assume that $0<H<\frac{1}{2q}$. Then, as $n\rightarrow \infty $,
it holds
\begin{equation}
2^{nqH-n}\,V_{n}^{(q)}(f)\,\,\,{\overset{L^{2}}{\longrightarrow }}\,\,\,%
\frac{(-1)^{q}}{2^{q}q!}\,\int_{0}^{1}f^{(q)}(B_{s})ds.  \label{t1}
\end{equation}

\item[2.]  Assume that $\frac{1}{2q}<H<1-\frac{1}{2q}$. Then, as $%
n\rightarrow \infty $, it holds
%%% debut
\begin{equation}
\big(B,2^{-n/2 }\,V_{n}^{(q)}(f)\big) \\ \,\,{\overset{\mathrm{Law}}{\longrightarrow }}
\,\,\,
\big(B,\sigma_{H,q}\int_{0}^{1}f(B_{s})dW_{s}\big),  \label{t2}
\end{equation}
where $W$ is a standard Brownian motion independent of $B$ and
\begin{equation}\label{hsmall}
\sigma_{H,q}=\sqrt{\frac{1}{2^qq!}\sum_{r\in\mathbb{Z}}\big(|r+1|^{2H}+|r-1|^{2H}-2|r|^{2H}\big)^q}.
\end{equation}
%%% fin

\item[3.] Assume that $H=1-\frac{1}{2q}$. Then, as $n\rightarrow \infty $,
it holds
\begin{equation}
\big(B,\frac{1}{\sqrt{n}}2^{-n/2}\,V_{n}^{(q)}(f)\big)\,\,\,{\overset{\mathrm{Law}}{%
\longrightarrow }}\,\,\,\big(B,\sigma_{1-1/(2q),q}\int_{0}^{1}f(B_{s})dW_{s}\big),
\end{equation}%
where $W$ is a standard Brownian motion independent of $B$ and
%%% debut
\begin{equation}\label{hcrit}
\sigma_{1-1/(2q),q}=\frac{2\log 2}{q!}\big(1-\frac1{2q}\big)^q\big(1-\frac1q\big)^q.
\end{equation}
%%% fin

\item[4.]  Assume that $H>1-\frac{1}{2q}$. Then, as $n\to\infty$, it holds
\begin{equation}
2^{nq(1-H)-n}\,V_{n}^{(q)}(f)\,\,\,{\overset{L^{2}}{\longrightarrow }}%
\,\,\,\int_{0}^{1}f(B_{s})dZ_{s}^{(q)},  \label{t3}
\end{equation}%
where $Z^{(q)}$ denotes the Hermite process
of order $q$ introduced in Definition \ref{d1} below.
\end{description}
\end{theorem}

\noindent \textbf{Remark 1.} When $q=1$, we have $V_{n}^{(1)}(f)=2^{-nH}
\sum_{k=1}^{2^{n}}f\left( B_{(k-1)2^{-n}}\right) \Delta B_{k2^{-n}}$.
For $H=\frac{1}{2}$, $2^{nH}V_n^{(1)}(f)$ converges in $L^{2}$ to the
It\^{o} stochastic integral $\int_{0}^{1}f(B_{s})dB_{s}$. For $H>\frac{1}{%
2}$, $2^{nH}V_n^{(1)}(f)$ converges in $L^{2}$ and almost surely  to the
Young integral $\int_{0}^{1}f(B_{s})dB_{s}$.
For $H<\frac12$, $2^{3nH-n}V_n^{(1)}(f)$ converges in $L^2$ to
$-\frac{1}{2}%
\int_0^1 f^{\prime}(B_s)ds$.

\vskip0.5cm

\noindent \textbf{Remark 2.} In the
critical case $H=\frac{1}{2q}$ ($q\geq 2$),
we conjecture the following asymptotic behavior: as $n\to\infty$,
\begin{equation}\label{conjecture}
\big(B,2^{-n/2}\,V_{n}^{(q)}(f)\big)\,\,\,{\overset{\mathrm{Law}}{\longrightarrow }}%
\,\,\,
\big(B,\sigma_{1/(2q),q}\int_{0}^{1}f(B_{s})dW_{s}+
\frac{(-1)^{q}}{2^{q}q!}\,\int_{0}^{1}f^{(q)}(B_{s})ds\big),
\end{equation}%
for $W$ a standard Brownian motion independent of $B$ and $\sigma_{1/(2q),q}$
%%% debut
the constant defined by (\ref{hsmall}).
%%% fin
Actually, (\ref{conjecture}) for $q=2$ and
$H=\frac14$ has been proved in
 \cite{burdzy-swanson,nourdin-nualart,nourdin-reveillac} after that
the first draft of the current paper have been submitted.
The reader is also referred to \cite{NoPe}
for the study of the weighted variations associated with iterated
Brownian motion, which is a non-Gaussian self-similar process of
order $\frac14$.

\vskip0.5cm

%%% debut
When $H$ is between $\frac14$ and $\frac34$, one can refine point {\bf 2} of Theorem \ref{main} as follows:
\begin{proposition}\label{joint}
Let $q\geq 2$ be an integer, $f:\mathbb{R}%
\rightarrow \mathbb{R}$ be a function such that $(\mathbf{H}_{q})$ holds
and assume that $H\in(\frac14,\frac34)$.
Then
\begin{eqnarray}\label{t2joint}
&&\left(B, 2^{-n/2 } V^{(2)}_n(f), \dots, 2^{-n/2 }\,V_{n}^{(q)}(f) \right)\\
&&\hskip3cm
{\overset{\mathrm{Law}}{\longrightarrow }}%
\,\,\,  \left( B, \sigma_{H,2}\int_{0}^{1}f(B_{s})dW^{(2)}_{s},\ldots,\sigma_{H,q} \int_{0}^{1}f(B_{s})dW^{(q)}_{s}\right),
\notag
\end{eqnarray}
where $(W^{(2)}, \dots, W^{(q)})$ is a $(q-1)$-dimensional  standard Brownian motion independent of $B$
and the $\sigma_{H,p}$'s, $2\le p\le q$, are given by (\ref{hsmall}).
\end{proposition}
%%% fin

\vskip0.5cm

 Theorem \ref{main} together with Proposition \ref{joint} allow to complete the missing
cases in the understanding of the asymptotic behavior of weighted {\it power}
variations of fractional Brownian motion:

\begin{corollary}
\label{complete-picture} Let $q\geq 2$ be an integer, and $f:\mathbb{R}%
\rightarrow \mathbb{R}$ be a function such that $(\mathbf{H}_{q})$ holds.
Then, as $n\rightarrow \infty $:

\begin{enumerate}
\item When $H>\frac{1}{2}$ and $q$ is odd,
\begin{equation}
2^{-nH}\sum_{k=1}^{2^{n}}f(B_{(k-1)2^{-n}})(2^{nH}\Delta B_{k2^{-n}})^{q}%
\overset{L^{2}}{\longrightarrow }q\mu _{q-1}\int_{0}^{1}f(B_{s})dB_{s}=q\mu
_{q-1}\int_{0}^{B_{1}}f(x)dx.  \label{c1}
\end{equation}

\item When $H<\frac{1}{4}$ and $q$ is even,
\begin{equation}
2^{2nH-n}\sum_{k=1}^{2^{n}}f(B_{(k-1)2^{-n}})\big[(2^{nH}\Delta
B_{k2^{-n}})^{q}-\mu _{q}\big]\overset{L^{2}}{\longrightarrow }
\frac{1}{4}%
\binom{q}{2}\mu _{q-2}\int_{0}^{1}f^{\prime \prime }(B_{s})ds.
  \label{c2}
\end{equation}%
(We recover (\ref{cv-L2}) by choosing $q=2$).

\item When $H= \frac{1}{4} $ and $q$ is even,
\begin{eqnarray}
&& \left(B,2^{-n/2}\sum_{k=1}^{2^{n}}f(B_{(k-1)2^{-n}})\big[(2^{n/4}\Delta
B_{k2^{-n}})^{q}-\mu _{q}\big]\right)
\overset{\mathrm{Law}}{\longrightarrow }  \notag
%%% debut
\left(B,\frac{1}{4}
\binom{q}{2}\mu _{q-2}\int_{0}^{1}f^{\prime \prime }(B_{s})ds\right.\\
&&\left.\hskip8cm+ \widetilde{\sigma}_{1/4,q} \int_0^1 f(B_s) dW_s\right),  \label{c3a1/4}
\end{eqnarray}
where $W$ is a standard Brownian motion independent of $B$
and $\widetilde{\sigma}_{1/4,q}$ is the constant given by (\ref{sity}) just below.
%%% fin

\item When $\frac{1}{4}<H<\frac{3}{4}$ and $q$ is even,
\begin{equation}
\left(B,2^{-n/2}\sum_{k=1}^{2^{n}}f(B_{(k-1)2^{-n}})\big[(2^{nH}\Delta
B_{k2^{-n}})^{q}-\mu _{q}\big]\right)\overset{\mathrm{Law}}{\longrightarrow }\left(B,\widetilde{\sigma}
_{H,q}\int_{0}^{1}f(B_{s})dW_{s}\right),  \label{c3}
\end{equation}%
%%% debut
for $W$ a standard Brownian motion independent of $B$
and
\begin{equation}\label{sity}
\widetilde{\sigma}_{H,q}=\sqrt{\sum_{p=2}^q p!\binom{q}{p}^2\mu_{q-p}^2\,2^{-p}\sum_{r\in\mathbb{Z}}
\big(|r+1|^{2H}+|r-1|^{2H}-2|r|^{2H}\big)^p}.
\end{equation}
%%% fin

\item When $H=\frac{3}{4}$ and $q$ is even,
\begin{equation}
\left(B,\frac{1}{\sqrt{n}}2^{-n/2}\sum_{k=1}^{2^{n}}f(B_{(k-1)2^{-n}})\big[(2^{nH}\Delta
B_{k2^{-n}})^{q}-\mu _{q}\big]\right)
\overset{\mathrm{Law}}{\longrightarrow }
\left(B,\widetilde{\sigma}
_{\frac34,q}\int_{0}^{1}f(B_{s})dW_{s}\right),  \label{c3a}
\end{equation}
%%% debut
for $W$ a standard Brownian motion
independent of $B$ and
$$
\widetilde{\sigma}_{\frac34,q}=\sqrt{\sum_{p=2}^q 2\log 2\,p!\binom{q}{p}^2\mu_{q-p}^2\big(
1-\frac1{2q}\big)^q\big(1-\frac{1}{q}\big)^q}.
$$
%%% fin

\item When $H>\frac{3}{4}$ and $q$ is even,
\begin{equation}
2^{n-2Hn}\sum_{k=1}^{2^{n}}f(B_{(k-1)2^{-n}})\big[(2^{nH}\Delta
B_{k2^{-n}})^{q}-\mu _{q}\big]\overset{L^{2}}{\longrightarrow }%
2\mu_{q-2}\binom{q}{2}
\int_{0}^{1}f(B_{s})dZ_{s}^{(2)},  \label{c4}
\end{equation}%
for $Z^{(2)}$ the
Hermite process introduced in Definition \ref{d1}.
\end{enumerate}
\end{corollary}

Finally, we can also give a new proof of the following result, stated and
proved by Gradinaru \textit{et al.} \cite{GNRV} and Cheridito
and Nualart \cite{CN} in a \textit{continuous} setting:

\begin{theorem}
\label{ito-thm} Assume that $H>\frac16$, and that $f:\mathbb{R}\rightarrow%
\mathbb{R}$ verifies \textrm{($\mathbf{H}_{\mathbf{6}}$)}. Then the limit in
probability, as $n\rightarrow \infty$, of the symmetric Riemann sums
\begin{equation}  \label{ito}
\frac12\sum_{k=1}^{2^n} \big(f^{\prime}(B_{k2^{-n}})+f^{%
\prime}(B_{(k-1)2^{-n}})\big) \,\Delta B_{k2^{-n}}
\end{equation}
exists and is given by $f(B_1)-f(0)$.
\end{theorem}

\noindent \textbf{Remark 3} When $H\le \frac16$, quantity
(\ref{ito}) does not converge
%%% debut
in probability
%%% fin
in general. As a counterexample, one
can consider the case where $f(x)=x^3$, see Gradinaru \textit{et
al.} \cite{GNRV} or Cheridito and Nualart \cite{CN}.

\section{Preliminaries and notation}

\label{sec2}

We briefly recall some basic facts about stochastic calculus with
respect to a fractional Brownian motion. One refers to \cite{Ncours}
for further details. Let $B=(B_{t})_{t\in \lbrack 0,1]}$ be a
fractional Brownian motion with Hurst parameter $H\in (0,1)$. That is, $B$
is a zero mean Gaussian process, defined on a complete probability
space $(\Omega,\mathcal{A},P)$, with the covariance function
$$
R_{H}(t,s)=E(B_tB_s)=\frac12\big(s^{2H}+t^{2H}-|t-s|^{2H}\big),\quad s,t\in[0,1].
$$
We suppose that $\mathcal{A}$ is the sigma-field generated by $B$.
Let $\mathscr{E}$ be the set of step functions on $[0,T]$, and
$\EuFrak H$ be the Hilbert space defined as the closure of
$\mathscr{E}$ with respect to the inner product
\begin{equation*}
\langle \mathbf{1}_{[0,t]},\mathbf{1}_{[0,s]}\rangle _{\EuFrak
H}=R_{H}(t,s).
\end{equation*}%
The mapping $\mathbf{1}_{[0,t]}\mapsto B_{t}$ can be extended to an
isometry between $\EuFrak H$ and the Gaussian space
$\mathcal{H}_{1}$ associated with $B$. We will denote this isometry
by $\varphi \mapsto B(\varphi )$.

Let $\mathscr{S}$ be the set of all smooth cylindrical random variables,
\textit{i.e.} of the form
\begin{equation*}
F=\phi (B_{t_{1}},\ldots ,B_{t_{m}})
\end{equation*}%
where $m\geq 1$, $\phi :\mathbb{R}^{m}\rightarrow \mathbb{R}\in \mathscr{C}%
_{b}^{\infty }$ and $0\leq t_{1}<\ldots <t_{m}\leq 1$. The derivative of $F$
with respect to $B$ is the element of $L^{2}(\Omega ,\EuFrak{H})$ defined by
\begin{equation*}
D_{s}F\;=\;\sum_{i=1}^{m}\frac{\partial \phi }{\partial x_{i}}%
(B_{t_{1}},\ldots ,B_{t_{m}})\mathbf{1}_{[0,t_{i}]}(s),\quad s\in \lbrack
0,1].
\end{equation*}%
In particular $D_{s}B_{t}=\mathbf{1}_{[0,t]}(s)$. For any integer $k\geq 1$,
we denote by ${\mathbb{D}}^{k,2}$ the closure of the set of smooth random
variables with respect to the norm
\begin{equation*}
\Vert F\Vert _{k,2}^{2}\;=\;E(F^{2}) +\sum_{j=1}^{k}E\left[
\|D^{j}F\|_{\EuFrak H^{\otimes j}}^{2}\right] .
\end{equation*}%
The Malliavin derivative $D$ satisfies the chain rule. If $\varphi :\mathbb{R}%
^{n}\rightarrow \mathbb{R}$ is $\mathscr{C}_{b}^{1}$ and if $%
(F_{i})_{i=1,\ldots ,n}$ is a sequence of elements of ${\mathbb{D}}^{1,2}$,
then $\varphi (F_{1},\ldots ,F_{n})\in {\mathbb{D}}^{1,2}$ and we have
\begin{equation*}
D\,\varphi (F_{1},\ldots ,F_{n})=\sum_{i=1}^{n}\frac{\partial \varphi }{%
\partial x_{i}}(F_{1},\ldots ,F_{n})DF_{i}.
\end{equation*}
We also have the following formula, which can easily be proved
by induction on $q$. Let $\varphi,\psi\in\mathscr{C}^q_b$ ($q\geq 1$),
and fix $0\leq u<v\leq 1$ and $0\leq s<t\leq 1$. Then $\varphi(B_t-B_s)\psi(B_v-B_u)\in\mathbb{D}^{q,2}$ and
\begin{equation}\label{der-multiple}
D^q\big(\varphi(B_t-B_s)\psi(B_v-B_u)\big)=\sum_{a=0}^{q}\binom{q}{a}\varphi^{(a)}(B_t-B_s)\psi^{(q-a)}(B_v-B_u)
{\bf 1}_{[s,t]}^{\otimes a}\widetilde{\otimes}{\bf 1}_{[u,v]}^{\otimes (q-a)},
\end{equation}
where $\widetilde{\otimes}$ means the symmetric tensor product.

The divergence operator $I$ is the adjoint of the derivative operator $%
D$. If a random variable $u\in L^{2}(\Omega ,\EuFrak H)$ belongs to the
domain of the divergence operator, that is, if it satisfies
\begin{equation*}
|E\langle DF,u\rangle _{\EuFrak H}|\leq c_{u}\,\sqrt{E(F^2)}\quad %
\mbox{for any }F\in {\mathscr{S}},
\end{equation*}%
then $I(u)$ is defined by the duality relationship
\begin{equation*}
E\big(FI(u)\big)=E\big(\langle DF,u\rangle _{\EuFrak H}\big),
\end{equation*}%
for every $F\in {\mathbb{D}}^{1,2}$.

For every $n\geq 1$, let $\mathcal{H}_{n}$ be the $n$th Wiener chaos
of $B,$
that is, the closed linear subspace of $L^{2}\left( \Omega ,\mathcal{A}%
,P\right) $ generated by the random variables $\{H_{n}\left( B\left(
h\right) \right) ,h\in \EuFrak H,\| h\| _{\EuFrak H}=1\}$, where $H_{n}$ is
the $n$th Hermite polynomial. The mapping $I_{n}(h^{\otimes
n})=n!H_{n}\left( B\left( h\right) \right) $ provides a linear isometry
between the symmetric tensor product $\EuFrak H^{\odot n}$ (equipped with the modified norm
$\|\cdot\|_{\EuFrak H^{\odot n}}=\frac{1}{\sqrt{n!}}
\|\cdot\|_{\EuFrak H^{\otimes n}}$) and $\mathcal{H}_{n}$.
For $H=\frac{1}{2}$, $I_{n}$ coincides with the multiple Wiener-It\^o
integral of order $n$. The following duality formula holds%
\begin{equation}
E\left( FI_{n}(h)\right) =E\left( \left\langle D^{n}F,h\right\rangle _{%
\EuFrak H^{\otimes n}}\right) ,  \label{dual}
\end{equation}%
for any element $h\in \EuFrak H^{\odot n}$ and any random variable $F\in
\mathbb{D}^{n,2}$.

Let $\{e_{k},\,k\geq 1\}$ be a complete orthonormal system in $\EuFrak H$.
Given $f\in \EuFrak H^{\odot n}$ and $g\in \EuFrak H^{\odot m}$, for every
$r=0,\ldots ,n\wedge m$, the contraction of $f$ and $g$ of order $r$
is the element of $\EuFrak H^{\otimes (n+m-2r)}$ defined by
$$
f\otimes _{r}g=\sum_{k_{1},\ldots ,k_{r}=1}^{\infty }\langle
f,e_{k_{1}}\otimes \ldots \otimes e_{k_{r}}\rangle _{\EuFrak H^{\otimes
r}}\otimes \langle g,e_{k_{1}}\otimes \ldots \otimes e_{k_{r}}
\rangle_{\EuFrak H^{\otimes r}}.
$$
Notice that $f\otimes _{r}g$ is not necessarily symmetric: we denote its
symmetrization by $f\widetilde{\otimes }_{r}g\in \EuFrak H^{\odot (n+m-2r)}$.
We have the following product formula: if $f\in \EuFrak
H^{\odot n}$ and $g\in \EuFrak H^{\odot m}$ then
\begin{eqnarray}\label{multiplication}
I_n(f) I_m(g) = \sum_{r=0}^{n \wedge m} r! {n \choose r}{ m \choose r} I_{n+m-2r} (f\widetilde{\otimes}_{r}g).
\end{eqnarray}

We recall the following simple formula for any $s<t$ and $u<v$:%
\begin{equation}
E\left( (B_{t}-B_{s})(B_{v}-B_{u})\right) =\frac{1}{2}\left(
|t-v|^{2H}+|s-u|^{2H}-|t-u|^{2H}-|s-v|^{2H}\right) .  \label{eq1}
\end{equation}%
We will also need the following lemmas:

\begin{lemma}
\label{lemma1}
\begin{enumerate}
\item Let $s<t$ belong to $[0,1]$. Then, if $H<1/2$, one has
\begin{equation}
\left| E\big(B_{u}(B_{t}-B_{s})\big)\right| \leq (t-s)^{2H}  \label{aux1}
\end{equation}
for all $u\in[0,1]$.

\item For  all $H\in (0,1)$,
\begin{equation}
\sum_{k,l =1}^{2^{n}}\left| E\left( B_{(k-1)2^{-n}}\,\Delta B_{l
2^{-n}}\right) \right| = O(2^{n}).  \label{aux2}
\end{equation}

\item For any $r\geq 1$, we have, if $H<1-\frac{1}{2r}$,
\begin{equation}
\sum_{k,l =1}^{2^{n}}\left| E\left( \Delta B_{k2^{-n}}\,\Delta B_{l
2^{-n}}\right) \right| ^{r}
=O(2^{n-2rHn}).  \label{aux4}
\end{equation}

\item For any $r\geq 1$, we have, if $H=1-\frac{1}{2r}$,
\begin{equation}
\sum_{k,l =1}^{2^{n}}\left| E\left( \Delta B_{k2^{-n}}\,\Delta B_{l
2^{-n}}\right) \right| ^{r}=O(n2^{2n-2rn}).  \label{aux6}
\end{equation}
\end{enumerate}
\end{lemma}
\textbf{Proof}: To prove inequality (\ref{aux1}), we just write
\begin{equation*}
E(B_{u}(B_{t}-B_{s}))=\frac{1}{2}(t^{2H}-s^{2H})+\frac{1}{2}\left(
|s-u|^{2H}-|t-u|^{2H}\right) ,
\end{equation*}%
and observe that we have $|b^{2H}-a^{2H}|\leq |b-a|^{2H}$ for any $a,b\in
\lbrack 0,1]$, because $H<\frac{1}{2}$. To show (\ref{aux2}) using (\ref{eq1}), we write
\begin{eqnarray*}
\sum_{k,l =1}^{2^{n}}\left| E\left( B_{(k-1)2^{-n}}\,\Delta B_{l
2^{-n}}\right) \right| &=&2^{-2Hn-1}\sum_{k,l =1}^{2^{n}}\left|
|l-1|^{2H}-l^{2H}-|l -k+1|^{2H}+|l -k|^{2H}\right| \\
&\leq &C2^{n},
\end{eqnarray*}
the last bound coming from a telescoping sum argument.
Finally, to show (\ref{aux4}) and (\ref{aux6}), we write
\begin{eqnarray*}
\sum_{k,l =1}^{2^{n}}\left| E\left( \Delta B_{k2^{-n}}\,\Delta B_{l
2^{-n}}\right) \right| ^{r}
&=&2^{-2nrH-r}\sum_{k,l =1}^{2^{n}}\big||k-l +1|^{2H}+|k-l
-1|^{2H}-2|k-l |^{2H}\big|^{r} \\
&\leq& 2^{n-2nrH-r}\sum_{p=-\infty }^{\infty }\big|%
|p+1|^{2H}+|p-1|^{2H}-2|p|^{2H}\big|^{r},
\end{eqnarray*}%
and observe that, since the function $\big|
|p+1|^{2H}+|p-1|^{2H}-2|p|^{2H}\big| $ behaves as $C_{H}p^{2H-2}$ for large $p$,  the series in the right-hand side is convergent because $%
H<1-\frac{1}{2r}$. In the critical case $H=1-\frac{1}{2r}$, this
series is divergent,
and%
\begin{equation*}
\sum_{p=-2^{n}}^{2^{n}}\big||p+1|^{2H}+|p-1|^{2H}-2|p|^{2H}\big|^{r}
\end{equation*}%
behaves as a constant time $n$. \qed

\vskip0.5cm

\begin{lemma}
\label{lemma2} Assume that $H>\frac12$.

\begin{enumerate}
\item Let $s<t$ belong to $[0,1]$. Then
\begin{equation}
\big|E\big(B_{u}(B_{t}-B_{s})\big)\big|\leq 2H(t-s)  \label{aux5}
\end{equation}
for all $u\in[0,1]$.

\item Assume that $H>1-\frac{1}{2l}$ for some $l\geq 1$. Let $u<v$ and $s<t$
belong to $[0,1]$. Then
\begin{equation}
\left| E(B_{u}-B_{v})(B_{t}-B_{s})\right| \leq H(2H-1)\left( \frac{2}{%
2Hl+1-2l}\right) ^{\frac{1}{l}}(u-v)^{\frac{l-1}{l}}(t-s).
\label{boundHgrand}
\end{equation}

\item Assume that $H>1-\frac{1}{2l}$ for some $l\ge 1$. Then
\begin{equation}
\sum_{i,j=1}^{2^n} \left| E\big(\Delta B_{i2^{-n}}\,\Delta B_{j2^{-n}}\big) %
\right|^l = O(2^{2n-2ln}).  \label{aux4new}
\end{equation}
\end{enumerate}
\end{lemma}
\textbf{Proof: } We have
\begin{equation*}
E\big(B_u(B_t-B_s)\big)=\frac12\big(t^{2H}-s^{2H}\big)+\frac12\big(%
|s-u|^{2H}-|t-u|^{2H}\big).
\end{equation*}
But, when $0\le a<b\le 1$:
\begin{equation*}
b^{2H}-a^{2H}=2H\int_0^{b-a} (u+a)^{2H-1}du \le 2H\,b^{2H-1}\,(b-a)\le
2H(b-a).
\end{equation*}
Thus, $|b^{2H}-a^{2H}|\le 2H|b-a|$ and the first point follows.

\vskip0.3cm

Concerning the second point, using H\"older inequality, we can write
\begin{eqnarray*}
\left| E(B_{u}-B_{v})(B_{t}-B_{s}) \right|&=&H(2H-1)
\int_{u}^{v}\int_{s}^{t} \vert y-x\vert ^{2H-2} dydx \\
&\leq & H(2H-1) \vert u-v\vert ^{\frac{l-1}{l}}\left( \int_{0}^{1} \left(
\int_{s}^{t} \vert y-x\vert ^{2H-2} dy \right) ^{l}dx\right) ^{\frac{1}{l} }
\\
&\leq & H(2H-1) \vert u-v\vert ^{\frac{l-1}{l}}\vert t-s\vert ^{\frac{l-1}{l}%
}\left( \int_{0}^{1}\int_{s}^{t}\vert y-x\vert ^{(2H-2)l} dy dx\right) ^{%
\frac{1}{l}}.
\end{eqnarray*}
Denote by $H^{\prime}= 1+ (H-1)l$ and observe that $H^{\prime}>\frac{1}{2} $
(because $H>1-\frac{1}{2l}$). Since $2H'-2=(2H-2)l$, we can write
\begin{equation*}
H^{\prime}(2H^{\prime}-1)\int_{0}^{1}\int_{s}^{t}\vert y-x\vert ^{(2H-2)l}
dy dx=E\left| B_{1}^{H^{\prime}} (B^{H^{\prime}}_{t}- B^{H^{\prime}}_{s}
)\right|\leq 2H^{\prime}\vert t-s \vert
\end{equation*}
by the first point of this lemma. This gives the desired bound.

\vskip0.3cm We prove now the third point. We have
\begin{eqnarray*}
\sum_{i,j=1}^{2^{n}}\left| E\big(\Delta B_{i2^{-n}}\,\Delta B_{j2^{-n}}\big)%
\right| ^{l} &=&2^{-2Hnl-l}\sum_{i,j=1}^{2^{n}}\left|
|i-j+1|^{2H}+|i-j-1|^{2H}-2|i-j|^{2H}\right| ^{l} \\
&\leq &2^{n-2Hnl+1-l}\sum_{k=-2^{n}+1}^{2^{n}-1}\left|
|k+1|^{2H}+|k-1|^{2H}-2|k|^{2H}\right| ^{l}
\end{eqnarray*}%
and the function $|k+1|^{2H}+|k-1|^{2H}-2|k|^{2H}$ behaves as $|k|^{2H-2}$
for large $k$. As a consequence, since $H>1-\frac{1}{2l}$, the sum
\begin{equation*}
\sum_{k=-2^n+1 }^{2^n-1 }\left| |k+1|^{2H}+|k-1|^{2H}-2|k|^{2H}\right| ^{l}
\end{equation*}%
behaves as $2^{(2H-2)ln+n}$ and the third point follows. \hfill %
\vrule width.25cm height.25cm depth0cm\smallskip

\vskip0.3cm

Now, let us introduce the Hermite process of order $q\geq 2$ appearing in (\ref{t3}). Fix $H>1/2$ and $t\in[0,1]$.
The sequence $\big(\varphi_n(t)\big)_{n\geq 1}$, defined as
\begin{equation*}
\varphi _{n}(t)=2^{nq-n}\frac1{q!}\sum_{j=1}^{\left[ 2^{n}t\right] }{\bf 1}_{[(j-1)2^{-n},j2^{-n}]}^{\otimes q},
\end{equation*}
is a Cauchy sequence in the space $\EuFrak H^{\otimes q}$. Indeed, since $H>1/2$, we have
$$
\langle {\bf 1}_{[a,b]},{\bf 1}_{[u,v]}\rangle_{\EuFrak H}=E\big((B_b-B_a)(B_v-B_u)\big)=H(2H-1)\int_a^b\int_u^v|s-s'|^{2H-2}dsds',
$$
so that, for any $m\geq n$
\begin{equation*}
\left\langle \varphi _{n}(t),\varphi _{m}(t)\right\rangle _{\EuFrak H^{\otimes
q}}=\frac{H^q(2H-1)^{q}}{q!^2}2^{nq+mq-n-m}\sum_{j=1}^{\left[ 2^{m}t\right] }\sum_{k=1}^{\left[
2^{n}t\right] }\left(
\int_{(j-1)2^{-m}}^{j2^{-m}}\int_{(k-1)2^{-n}}^{k2^{-n}}|s-s^{\prime
}|^{2H-2}dsds^{\prime }\right) ^{q}.
\end{equation*}%
Hence
$$
\lim_{m,n\rightarrow \infty }\left\langle \varphi_n(t),
\varphi _{m}(t)\right\rangle _{\EuFrak H^{\otimes q}} =\frac{H^q(2H-1)^{q}}{q!^2}
 \int_{0}^{t}\int_{0}^{t}|s-s^{\prime }|^{(2H-2)q}dsds^{\prime }
=c_{q,H}t^{(2H-2)q+2},
$$
where $c_{q,H}=\frac{H^q(2H-1)^{q}}{q!^2(Hq-q+1)(2Hq-2q+1) }$.
Let us denote by $\mu^{(q)}_{t}$
%the uniform measure on the diagonal of $[0,t]^{q}$, normalized to have a
%total mass equal to $t/q!$.
%Then, 
the limit in $\EuFrak H^{\otimes q}$ of the sequence
of functions $\varphi _{n}(t)$.
% can be identified with $\mu^{(q)}_{t}$. Indeed, for 
For any $f\in\EuFrak H^{\otimes q}$, we have
\begin{eqnarray*}
&&\langle \varphi_n(t),f\rangle_{\EuFrak H^{\otimes q}}
=2^{nq-n}\frac1{q!}\sum_{j=1}^{[2^nt]}\langle {\bf 1}_{[(j-1)2^{-n},j2^{-n}]}^{\otimes q},f\rangle_{\EuFrak H^{\otimes q}}\\
&=&2^{nq-n}\frac1{q!}H^q(2H-1)^q\sum_{j=1}^{[2^nt]}\int_0^1 ds_1\int_{(j-1)2^{-n}}^{j2^{-n}}ds'_1|s_1-s'_1|^{2H-2}\ldots\\
&&\hskip3cm\times
\int_0^1 ds_q \int_{(j-1)2^{-n}}^{j2^{-n}}ds'_q|s_q-s'_q|^{2H-2}f(s_1,\ldots,s_q)\\
&\underset{n\to\infty}{\longrightarrow}&
\frac1{q!}H^q(2H-1)^q \int_0^t ds'
\int_{[0,1]^q} ds_1\ldots ds_q|s_1-s'|^{2H-2}\ldots
|s_q-s'|^{2H-2}f(s_1,\ldots,s_q)\\
&&=\langle\mu_t^{(q)},f\rangle_{\EuFrak H^{\otimes q}}.
\end{eqnarray*}

\begin{definition}
\label{d1} Fix $q\geq 2$ and $H>1/2$. The \textit{Hermite
process} $Z^{(q)}=(Z^{(q)}_t)_{t\in[0,1]}$ of order $q$ is defined by $Z_{t}^{(q)}=I_{q}(\mu^{(q)}_{t})$
for $t\in[0,1]$.
\end{definition}

Let $Z^{(q)}_{n}$
be the process defined by $Z^{(q)}_n(t)=I_q(\varphi_n(t))$ for $t\in[0,1]$.
By construction, it is clear that $Z_n^{(q)}(t)\overset{L^2}{\longrightarrow} Z^{(q)}(t)$ as $n\to\infty$,
for all fixed $t\in[0,1]$. On the other hand,
it follows, from Taqqu
\cite{Taqqu}
and Dobrushin and Major \cite{DM}, that $Z^{(q)}_{n}$
converges in law to the ``standard'' and historical $q$th Hermite process,
defined through its moving average representation as a multiple integral
with respect to a Wiener process with time horizon $\mathbb{R}$.
In particular, the process introduced in Definition \ref{d1}
has the same finite dimensional distributions as the historical Hermite
process.

Let us finally mention that
it can be easily seen that $Z^{(q)}$ is $%
q(H-1)+1$ self-similar, has stationary increments and admits moments
of all orders. Moreover, it has H\"{o}lder continuous paths
of order strictly less than $q(H-1)+1$.
For further results, we refer to Tudor \cite{cip-rosenblatt}.

\section{Proof of the main results}

In this section we will provide the proofs of the main results. For
notational convenience, from now on, we write $\varepsilon _{(k-1)2^{-n}}$ (resp. $\delta _{k2^{-n}}$%
) instead of $\mathbf{1}_{[0,(k-1)2^{-n}]}$ (resp. $\mathbf{1}%
_{[(k-1)2^{-n},k2^{-n}]}$). \ The following proposition provides information
on the asymptotic behavior of \ $E\left( V_{n}^{(q)}(f)^{2}\right) $, as $n$
tends to infinity, \ for $H\leq 1-\frac{1}{2q}$.

\begin{proposition}
\label{p1}Fix an integer $q\geq 2$.  Suppose that $f$ satisfies  ($\mathbf{H}_q$). Then, if $H\leq\frac{1}{2q}$, then%
\begin{equation}
E\left( V_{n}^{(q)}(f)^{2}\right) = O(2^{n(-2Hq+2)}).  \label{r1}
\end{equation}%
If $\frac{1}{2q}\leq H<1-\frac{1}{2q}$, then%
\begin{equation}
E\left( V_{n}^{(q)}(f)^{2}\right) = O(2^{n}).  \label{r1a}
\end{equation}%
Finally, if $H=1-\frac{1}{2q}$, then%
\begin{equation}
E\left( V_{n}^{(q)}(f)^{2}\right) = O(n2^{n}).  \label{r1b}
\end{equation}
\end{proposition}
\begin{proof}
Using the relation between Hermite polynomials and multiple stochastic
integrals, we have $H_{q}\left( 2^{nH}\Delta B_{k2^{-n}}\right) =\frac{1}{q!}%
2^{qnH}I_{q}\left( \delta _{k2^{-n}}^{\otimes q}\,\right) $. In this way we
obtain
\begin{eqnarray*}
&&E\left( V_{n}^{(q)}(f)^{2}\right) \\
&=&\sum_{k,l =1}^{2^{n}}E\left\{ f(B_{(k-1)2^{-n}})\,f(B_{(l
-1)2^{-n}})\,H_{q}\left( 2^{nH}\Delta B_{k2^{-n}}\right) \,H_{q}\left(
2^{nH}\Delta B_{l 2^{-n}}\right) \right\} \\
&=&\frac{1}{q!^{2}}\,2^{2Hqn}\sum_{k,l =1}^{2^{n}}E\left\{
f(B_{(k-1)2^{-n}})\,f(B_{(l -1)2^{-n}})\,I_{q}\left( \delta
_{k2^{-n}}^{\otimes q}\,\right) I_{q}\left( \delta _{l 2^{-n}}^{\otimes
q}\right) \right\} .
\end{eqnarray*}%
Now we apply the product formula (\ref{multiplication}) for multiple stochastic integrals
and the duality relationship (\ref{dual})
between the multiple stochastic integral $I_{N}$ and
the iterated derivative operator $D^{N}$, obtaining%
\begin{eqnarray*}
&&E\left( V_{n}^{(q)}(f)^{2}\right) \\
&=&\frac{2^{2Hqn}}{q!^2}\sum_{k,l =1}^{2^{n}}\sum_{r=0}^{q}r!\binom{q}{r%
}^{2} \\
&&\times E\left\{ f(B_{(k-1)2^{-n}})\,f(B_{(l
-1)2^{-n}})\,I_{2q-2r}\left( \delta _{k2^{-n}}^{\otimes q-r}\widetilde{%
\otimes }\delta _{l 2^{-n}}^{\otimes q-r}\right) \right\} \langle \delta
_{k2^{-n}},\delta _{l 2^{-n}}\rangle _{\EuFrak H}^{r} \\
&=&2^{2Hqn}\sum_{k,l =1}^{2^{n}}\sum_{r=0}^{q}\frac{1}{r!
(q-r)!^{2}} \\
&&\times E\left\{ \left\langle D^{2q-2r}\left(
f(B_{(k-1)2^{-n}})\,f(B_{(l -1)2^{-n}})\right) ,\delta
_{k2^{-n}}^{\otimes q-r}\widetilde{\otimes }\delta _{l 2^{-n}}^{\otimes
q-r}\right\rangle _{\EuFrak H^{\otimes (2q-2r)}}\right\} \langle \delta
_{k2^{-n}},\delta _{l 2^{-n}}\rangle _{\EuFrak H}^{r},
\end{eqnarray*}%
where $\widetilde{\otimes }$ denotes the symmetrization of the
tensor product. By (\ref{der-multiple}), the derivative of the product \ $D^{2q-2r}\left(
f(B_{(k-1)2^{-n}})\,f(B_{(l -1)2^{-n}})\right) $ is equal to
a sum
of derivatives:%
\begin{eqnarray*}
D^{2q-2r}\left( f(B_{(k-1)2^{-n}})\,f(B_{(l -1)2^{-n}})\right)
&=&\sum_{a+b=2q-2r}\,\,\,\,f^{(a)}(B_{(k-1)2^{-n}})\,f^{(b)}(B_{(l
-1)2^{-n}}) \\
&&\times \frac{(2q-2r)!}{a!b!}\left( \varepsilon _{(k-1)2^{-n}}^{\otimes a}\widetilde{\otimes} \varepsilon
_{(l -1)2^{-n}}^{\otimes b}\right) .
\end{eqnarray*}%
We make the decomposition%
\begin{equation}
E\left( V_{n}^{(q)}(f)^{2}\right) =A_{n}+B_{n}+C_{n}+D_{n},  \label{dec1}
\end{equation}%
where
\begin{equation*}
A_{n}=\frac{2^{2Hqn}}{q!^2}\sum_{k,l =1}^{2^{n}}E\left\{
f^{(q)}(B_{(k-1)2^{-n}})\,f^{(q)}(B_{(l -1)2^{-n}})\right\} \langle
\varepsilon _{(k-1)2^{-n}},\delta _{k2^{-n}}\rangle ^{q}\,\,\langle
\varepsilon _{(l -1)2^{-n}},\delta _{l 2^{-n}}\rangle ^{q},
\end{equation*}%
\begin{eqnarray*}
B_{n} &=&2^{2Hqn}\sum_{\substack{ c+d+e+f=2q  \\ d+e\geq 1}}\,\,\,\sum_{k,l
=1}^{2^{n}}E\left\{ f^{(q)}(B_{(k-1)2^{-n}})\,f^{(q)}(B_{(l
-1)2^{-n}})\right\} \alpha (c,d,e,f) \\
&&\times \langle \varepsilon _{(k-1)2^{-n}},\delta _{k2^{-n}}\rangle _{%
\EuFrak H}^{c}\langle \varepsilon _{(k-1)2^{-n}},\delta _{l
2^{-n}}\rangle _{\EuFrak H}^{d}\,\langle \varepsilon _{(l
-1)2^{-n}},\delta _{k2^{-n}}\rangle _{\EuFrak H}^{e}\,\langle \varepsilon
_{(l -1)2^{-n}},\delta _{l 2^{-n}}\rangle _{\EuFrak H}^{f},\\
C_{n}&=&2^{2Hqn}\sum_{\substack{ a+b=2q  \\ (a,b)\neq (q,q)}}%
\,\,\,\,\sum_{k,l =1}^{2^{n}}E\left\{
f^{(a)}(B_{(k-1)2^{-n}})\,f^{(b)}(B_{(l -1)2^{-n}})\right\} \,\frac{(2q)!}{%
q!^{2}a!b!} \\
&&\times \langle \varepsilon _{(k-1)2^{-n}}^{\otimes a}\widetilde{\otimes }%
\varepsilon _{(l -1)2^{-n}}^{\otimes b},\delta _{k2^{-n}}^{\otimes q}%
\widetilde{\otimes }\delta _{l 2^{-n}}^{\otimes q}\rangle _{\EuFrak %
H^{\otimes (2q)}},
\end{eqnarray*}%
and%
\begin{eqnarray*}
D_{n}&=&2^{2Hqn}\sum_{r=1}^{q}\sum_{a+b=2q-2r}\,\,\,\,\sum_{k,l
=1}^{2^{n}}E\left\{ f^{(a)}(B_{(k-1)2^{-n}})\,f^{(b)}(B_{(l
-1)2^{-n}})\right\} \,\frac{(2q-2r)!}{r!(q-r)!^{2}a!b!} \\
&&\times \langle \varepsilon _{(k-1)2^{-n}}^{\otimes a}\widetilde{\otimes }%
\varepsilon _{(l -1)2^{-n}}^{\otimes b},\delta _{k2^{-n}}^{\otimes q-r}%
\widetilde{\otimes }\delta _{l 2^{-n}}^{\otimes q-r}\rangle _{\EuFrak %
H^{\otimes (2q-2r)}}\left\langle \delta _{k2^{-n}},\delta _{l
2^{-n}}\right\rangle _{\EuFrak H}^{r},
\end{eqnarray*}%
for some combinatorial constants $\alpha (c,d,e,f)$. That is, \ $A_{n}$
and $B_{n}$ contain all the terms with $r=0$ and $(a,b)=(q,q)$; $C_{n}$
contains the terms with $r=0$ and \ $(a,b)\neq (q,q)$; and $D_{n}$ contains
the remaining terms.

For any integer $r\geq 1$, we set%
\begin{eqnarray}
\alpha _{n} &=&\sup_{k,l=1,\ldots,2^n}\left| \langle \varepsilon
_{(k-1)2^{-n}},\delta _{l 2^{-n}}\rangle _{\EuFrak H}\right|,
\label{alpha} \\
\beta _{r,n} &=&\sum_{k,l =1}^{2^{n}}\left| \left\langle \delta
_{k2^{-n}},\delta _{l 2^{-n}}\right\rangle _{\EuFrak H}\right|
_{\ }^{r}, \label{beta}\\
\gamma _{n} &=&\sum_{k,l =1}^{2^{n}}\left| \langle \varepsilon
_{(k-1)2^{-n}},\delta _{l 2^{-n}}\rangle _{\EuFrak H}\right|.
\end{eqnarray}%

Then, under assumption (${\bf H}_{\bf q}$), we have the following estimates:%
\begin{eqnarray*}
\left| A_{n}\right|  &\leq &C2^{2Hqn+2n}(\alpha_n)^{2q}, \\
\left| B_{n}\right| +\left| C_{n}\right|  &\leq &C2^{2Hqn}(\alpha_n)^{2q-1}\gamma _{n}, \\
\left| D_{n}\right|  &\leq &C2^{2Hqn}\sum_{r=1}^{q}(\alpha_{n})^{2q-2r}\beta_{r,n},
\end{eqnarray*}%
where $C$ is a constant depending only on $q$ and the function $f$. Notice that the
second inequality follows from the fact that when $(a,b)\neq (q,q)$, or
$(a,b)=(q,q)$ and $c+d+e+f=2q$ with $d\geq 1$ or $e\geq 1$, there will be at
least a factor of the form $\langle \varepsilon _{(k-1)2^{-n}},\delta _{l
2^{-n}}\rangle _{\EuFrak H}$ in the expression of $B_{n}$ or $C_{n}$.

In the case $H<\frac12$, we have  by (\ref{aux1}) that $\alpha _{n}\leq 2^{-2nH}$,
by (\ref{aux4}) that $\beta _{r,n}\leq C2^{n-2rHn}$, and by (\ref{aux2}) that
$\gamma _{n}\leq C2^{n}$. As a consequence, we obtain%
\begin{eqnarray}
\left| A_{n}\right| &\leq &C2^{n(-2Hq+2)},  \label{r2} \\
\left| B_{n}\right| +\left| C_{n}\right| &\leq &C2^{n(-2Hq+2H+1)},
\label{r3} \\
\left| D_{n}\right| &\leq &C\sum_{r=1}^{q}2^{n(-2(q-r)H+1)},  \label{r4}
\end{eqnarray}
which implies the estimates (\ref{r1}) and (\ref{r1a}).

\vskip0.3cm

In the case $\frac{1}{2}\leq H< 1-\frac{1}{2q}$, we have by (\ref{aux5}) that
$\alpha_{n}\leq C2^{-n}$, by (\ref{aux4}) that $\beta _{r,n}\leq C2^{n-2rHn}$, and by (\ref{aux2})
that $\gamma _{n}\leq C2^{n}$. As a consequence, we obtain
\begin{eqnarray*}
\left| A_{n}\right| +\left| B_{n}\right| +\left| C_{n}\right|&\leq &C2^{n(2q(H-1)+2)}, \\
\left| D_{n}\right| &\leq &C\sum_{r=1}^{q}2^{n((2q-2r)(H-1)+1)},
\end{eqnarray*}
which also implies (\ref{r1a}).

Finally, if  $H=1-\frac{1}{2q}$, we have by (\ref{aux5}) that $\alpha_{n}\leq C2^{-n}$,
by (\ref{aux6}) that $\beta_{r,n}\leq Cn2^{2n-2rn}$, and by (\ref{aux2}) that $\gamma _{n}\leq C2^{n}$.
As a consequence, we obtain
\begin{eqnarray*}
\left| A_{n}\right|+\left| B_{n}\right| +\left| C_{n}\right| &\leq &C2^{n}, \\
\left| D_{n}\right| &\leq &C\sum_{r=1}^{q}n2^{n\frac{r}q},
\end{eqnarray*}
which implies (\ref{r1b}).
\end{proof}

\subsection{Proof of Theorem \ref{main} in the case $0<H<\frac{1}{2q}$}

In this subsection we are going to prove the first point of Theorem \ref%
{main}. The proof will be done in three steps. Set
$V_{1,n}^{(q)}(f)=2^{n(qH-1)}V_{n}^{(q)}(f)$.
We first study the asymptotic behavior of $E\left(
V_{1,n}^{(q)}(f)^{2}\right) $, using Proposition \ref{p1}.

\textit{Step 1}. The decomposition (\ref{dec1}) leads to
\begin{equation*}
E\left( V_{1,n}^{(q)}(f)^{2}\right) =2^{2n(qH-1)}\left(
A_{n}+B_{n}+C_{n}+D_{n}\right) .
\end{equation*}
From the estimate (\ref{r3}) we obtain
$
2^{2n(qH-1)}\left( \left| B_{n}\right| +\left| C_{n}\right| \right) \leq
C2^{n(2H-1)},
$
which converges to zero as $n$ goes to infinity since $H<\frac{1}{2q}<
\frac{1}{2}$. On the other hand (\ref{r4}) yields
\begin{equation*}
2^{2n(qH-1)}\left| D_{n}\right| \leq C\sum_{r=1}^{q}2^{n(2rH-1)},
\end{equation*}
which tends to zero as $n$ goes to infinity since $2rH-1\leq 2qH-1<0$ for all $r=1,\ldots,q$.

In order to handle the term $A_{n}$, we make use of the following estimate,
which follows from (\ref{aux1}) and (\ref{eq1}):
\begin{eqnarray}
&&\left| \langle \varepsilon _{(k-1)2^{-n}},\delta _{k2^{-n}}\rangle _{%
\EuFrak H}^{q}-\left( -\frac{2^{-2Hn}}{2}\right) ^{q}\right|  \notag \\
&=&\left| \langle \varepsilon _{(k-1)2^{-n}},\delta _{k2^{-n}}\rangle _{%
\EuFrak H}+\frac{2^{-2Hn}}{2}\right| \,\left| \sum_{s=0}^{q-1}\langle
\varepsilon _{(k-1)2^{-n}},\delta _{k2^{-n}}\rangle _{\EuFrak H}^{s}\,\left(
-\frac{2^{-2Hn}}{2}\right) ^{q-1-s}\right|  \notag \\
&\leq &C\left( k^{2H}-(k-1)^{2H}\right) 2^{-2Hqn}.  \label{contr}
\end{eqnarray}%
Thus,
\begin{eqnarray*}
&&\left| \frac{2^{4Hqn-2n}}{q!^2}\sum_{k,l =1}^{2^{n}}E\left\{
f^{(q)}(B_{(k-1)2^{-n}})\,f^{(q)}(B_{(l -1)2^{-n}})\right\} \langle
\varepsilon _{(k-1)2^{-n}},\delta _{k2^{-n}}\rangle _{\EuFrak %
H}^{q}\,\,\langle \varepsilon _{(l -1)2^{-n}},\delta _{l
2^{-n}}\rangle _{\EuFrak H}^{q}\right. \\
&&\left. \hskip4cm-\,\frac{2^{-2n-2q}}{q!^2}\sum_{k,l =1}^{2^{n}}E\left\{
f^{(q)}(B_{(k-1)2^{-n}})\,f^{(q)}(B_{(l -1)2^{-n}})\right\} \right| \leq
C2^{2Hn-n},
\end{eqnarray*}%
which implies, as $n\rightarrow \infty $:
\begin{equation}
E\big(V_{1,n}^{(q)}(f)^{2}\big)=\frac{2^{-2n-2q}}{q!^2}\sum_{k,l =1}^{2^{n}}E\left\{
f^{(q)}(B_{(k-1)2^{-n}})\,f^{(q)}(B_{(l -1)2^{-n}})\right\} +o(1).
\label{A}
\end{equation}

\textit{Step 2: }We need the asymptotic behavior of the double product%
\begin{equation*}
J_{n}:=E\left( V_{1,n}^{(q)}(f)\times 2^{-n}\sum_{l =1}^{2^{n}}f^{(q)}(B_{(l
-1)2^{-n}})\right) .
\end{equation*}
Using the same arguments as in Step 1 we obtain
\begin{eqnarray*}
J_{n} &=&2^{Hqn-2n}\sum_{k,l =1}^{2^{n}}E\left\{
f(B_{(k-1)2^{-n}})\,f^{(q)}(B_{(l -1)2^{-n}})\,H_{q}\big(2^{nH}\Delta
B_{k2^{-n}}\big)\right\} \\
&=&\frac{1}{q!}2^{2Hqn-2n}\sum_{k,l =1}^{2^{n}}E\left\{
f(B_{(k-1)2^{-n}})\,f^{(q)}(B_{(l -1)2^{-n}})\,I_{q}\big(\delta
_{k2^{-n}}^{\otimes q}\big)\right\} \\
&=&\frac{1}{q!}2^{2Hqn-2n}\sum_{k,l =1}^{2^{n}}E\left\{ \left\langle D^{q}%
\big(f(B_{(k-1)2^{-n}})\,f^{(q)}(B_{(l -1)2^{-n}})\big),\delta
_{k2^{-n}}^{\otimes q}\right\rangle_{\HH^{\otimes q}} \right\} \\
&=&2^{2Hqn-2n}\sum_{k,l =1}^{2^{n}}\sum_{a=0}^{q}\frac{1}{a!(q-a)!}%
E\left\{ f^{(a)}(B_{(k-1)2^{-n}})\,f^{(2q-a)}(B_{(l -1)2^{-n}})\right\} \,
\\
&&\times \langle \varepsilon _{(k-1)2^{-n}},\delta _{k2^{-n}}\rangle _{%
\EuFrak H}^{a}\,\langle \varepsilon _{(l -1)2^{-n}},\delta
_{k2^{-n}}\rangle _{\EuFrak H}^{q-a}.
\end{eqnarray*}%
It turns out that only the term with $a=q$ will contribute to the limit as $%
n $ tends to infinity. For this reason we make the decomposition%
\begin{equation*}
J_{n}=2^{2Hqn-2n}\sum_{k,l =1}^{2^{n}}\frac{1}{q!}E\left\{
f^{(q)}(B_{(k-1)2^{-n}})\,f^{(q)}(B_{(l -1)2^{-n}})\right\} \,\langle
\varepsilon _{(k-1)2^{-n}},\delta _{k2^{-n}}\rangle _{\EuFrak H}^{q}+S_{n},
\end{equation*}%
where
\begin{eqnarray*}
&&S_{n}=2^{2Hqn-2n}\sum_{k,l =1}^{2^{n}}\langle \varepsilon _{(l
-1)2^{-n}},\delta _{k2^{-n}}\rangle _{\EuFrak H}\sum_{a=0}^{q-1}\frac{1}{a!(q-a)!}%
E\left\{ f^{(a)}(B_{(k-1)2^{-n}})\,f^{(2q-a)}(B_{(l
-1)2^{-n}})\right\} \\
&&\times \langle \varepsilon _{(k-1)2^{-n}},\delta _{k2^{-n}}\rangle _{%
\EuFrak H}^{a}\,\langle \varepsilon _{(l -1)2^{-n}},\delta
_{k2^{-n}}\rangle _{\EuFrak H}^{q-a-1}.
\end{eqnarray*}%
By (\ref{aux1}) and (\ref{aux2}), we have
\begin{equation*}
|S_{n}|\leq C2^{2Hn-n},
\end{equation*}%
which tends to zero as $n$ goes to infinity. Moreover, by (\ref{contr}), we
have
\begin{eqnarray*}
&&\left| \frac{2^{2Hqn-2n}}{q!}\sum_{k,l =1}^{2^{n}}E\left\{
f^{(q)}(B_{(k-1)2^{-n}})\,f^{(q)}(B_{(l -1)2^{-n}})\right\} \,\langle
\varepsilon _{(k-1)2^{-n}},\delta _{k2^{-n}}\rangle _{\EuFrak H}^{q}\right.
\\
&&\left. -(-1)^{q}\,\frac{2^{-2n-q}}{q!}\sum_{k,l =1}^{2^{n}}E\left\{
f^{(q)}(B_{(k-1)2^{-n}})\,f^{(q)}(B_{(l -1)2^{-n}})\right\} \right| \leq
C\,2^{2Hn-n},
\end{eqnarray*}%
which also tends to zero as $n$ goes to infinity. Thus, finally, as $%
n\rightarrow \infty $:
\begin{equation}
J_{n}=(-1)^{q}\,\frac{2^{-2n-q}}{q!}\sum_{k,l =1}^{2^{n}}E\left\{
f^{(q)}(B_{(k-1)2^{-n}})\,f^{(q)}(B_{(l -1)2^{-n}})\right\} +o(1).
\label{Ba}
\end{equation}

\textit{Step 3: } By combining (\ref{A}) and (\ref{Ba}), we obtain that
\begin{equation*}
E\left| V_{1,n}^{(q)}(f)-\frac{(-1)^{q}}{2^{q}q!}\,2^{-n}%
\sum_{k=1}^{2^{n}}f^{(q)}(B_{(k-1)2^{-n}})\right| ^{2}=o(1),
\end{equation*}%
as $n\rightarrow \infty $. Thus, the proof of the first point of Theorem \ref%
{main} is done using a Riemann sum argument. \qed

\vskip0.3cm

\subsection{Proof of Theorem \ref{main} in the case
$H>1-\frac{1}{2q}$: the weighted non-central limit theorem}

We prove here that the sequence $V_{3,n}(f)$, given by
\begin{equation*}
V_{3,n}^{(q)}(f)=2^{n(1-H)q-n}\,V_{n}^{(q)}(f)=2^{qn-n}\frac{1}{q!}%
\sum_{k=1}^{2^{n}}f\left( B_{(k-1)2^{-n}}\right) I_{q}\left( \delta
_{k2^{-n}}^{\otimes q}\right),
\end{equation*}%
converges in $L^{2}$ as $n\rightarrow \infty$ to the pathwise
integral $\int_{0}^{1}f(B_{s})dZ_{s}^{(q)}$ with respect to the Hermite
process of order $q$ introduced in Definition \ref{d1}.

Observe first that, by construction of $Z^{(q)}$ (precisely, see the discussion before Definition
\ref{d1} in Section \ref{sec2}),
the desired result is in order when the function $f$ is identically one. More precisely:

\begin{lemma}
\label{lemma3} For each fixed $t\in \lbrack 0,1]$, the sequence
$2^{qn-n}\frac{1}{q!}\sum_{k=1}^{\left[ 2^{n}t\right] }I_{q}\left(
\delta _{k2^{-n}}^{\otimes q}\right)$
converges in $L^{2}$ to the Hermite random variable $Z_{t}^{(q)}$.
\end{lemma}

Now, consider the case of a general function $f$. We fix two integers $m\geq n$, and
decompose the sequence \ $V_{3,m}^{(q)}(f)\ $as follows:
\begin{equation*}
V_{3,m}^{(q)}(f)=A^{\left( m,n\right) }+B^{\left( m,n\right) },
\end{equation*}%
where
\begin{equation*}
A^{(m,n)}=\ \frac{1}{q!}2^{m(q-1)}\sum_{j=1}^{2^{n}}f\left(
B_{(j-1)2^{-n}}\right) \sum_{i=(j-1)2^{m-n}+1}^{j2^{m-n}}I_{q}\left( \delta
_{i2^{-m}}^{\otimes q}\right) ,
\end{equation*}%
and
\begin{equation*}
B^{\left( m,n\right) }=\frac{1}{q!}\,2^{m(q-1)}\sum_{j=1}^{2^{n}}%
\sum_{i=(j-1)2^{m-n}+1}^{j2^{m-n}}\Delta _{i,j}^{m,n}f(B)\,I_{q}\left( \delta
_{i2^{-m}}^{\otimes q}\right) ,
\end{equation*}%
with the notation $\Delta _{i,j}^{m,n}f(B)=f(B_{\left( i-1\right)
2^{-m}})-f(B_{\left( j-1\right) 2^{-n}})$. We shall study $A^{\left(
m,n\right) }$ and $B^{\left( m,n\right) }$ separately.

\vskip0.5cm

\textit{Study of }$A^{\left( m,n\right) }$. When $n$ is fixed, Lemma \ref%
{lemma3} yields that the random vector
\begin{equation*}
\left( \frac{1}{q!}2^{m(q-1)}\sum_{i=\left( j-1\right)
2^{m-n}+1}^{j2^{m-n}}I_{q}\left( \delta _{i2^{-m}}^{\otimes q}\right)
;j=1,\ldots ,2^{n}\right)
\end{equation*}%
converges in $L^{2}$, as $m\rightarrow \infty $, to the vector
\begin{equation*}
\left( Z_{j2^{-n}}^{(q)}-Z_{(j-1)2^{-n}}^{(q)};j=1,\ldots ,2^{n}\right) .
\end{equation*}%
Then, as $m\rightarrow \infty $, $A^{\left( m,n\right) }\overset{L^{2}}{\rightarrow }A^{\left( \infty ,n\right) }$, where
\begin{equation*}
A^{(\infty ,n)}:=\sum_{j=1}^{2^{n}}f(B_{(j-1)2^{-n}})\left(
Z_{j2^{-n}}^{(q)}-Z_{(j-1)2^{-n}}^{(q)}\right) .
\end{equation*}%
Finally, we claim that when $n$ tends to infinity, $A^{(\infty ,n)}$
converges in $L^{2}$ to $\int_{0}^{1}f\left( B_{s}\right) dZ_{s}^{(q)}$.
Indeed, observe that the stochastic integral $\int_{0}^{1}f\left( B_{s}\right)
dZ_{s}^{(q)}$ is a pathwise Young integral. So, to get the
convergence in $L^{2}$ it suffices to show that the sequence $A^{(\infty
,n)} $ is bounded in $L^{p}$ for some $p\geq 2$.  The integral $%
\int_{0}^{1}f\left( B_{s}\right) dZ_{s}^{(q)}$ has moments of all orders,
because for all $p\geq 2$%
\begin{equation*}
E\left[ \sup_{0\leq s<t\leq 1}\left( \frac{\left|
Z_{t}^{(q)}-Z_{s}^{(q)}\right| }{|t-s|^{\gamma }}\right) ^{p}\right] <\infty
\end{equation*}%
and%
\begin{equation*}
E\left[ \sup_{0\leq s<t\leq 1}\left( \frac{\left| B_{t}-B_{s}\right| }{%
|t-s|^{\beta }}\right) ^{p}\right] <\infty ,
\end{equation*}%
if $\gamma <q(H-1)+1$ and $\beta <H$. On the other hand, Young's inequality
implies
\begin{equation*}
\left| A^{(\infty ,n)}-\int_{0}^{1}f\left( B_{s}\right) dZ_{s}^{(q)}\right|
\leq c_{\rho,\nu}\mathrm{Var}_{\rho}\big(f(B)\big)\mathrm{Var}_{\nu}\big(Z^{(q)}\big),
\end{equation*}%
where $\mathrm{Var}_{\rho}$ denotes the variation of order $\rho$, and with $\rho,\nu>1$
such that $\frac{1}{\rho}+\frac{1}{\nu}>1$. Choosing $\rho>%
\frac{1}{H}$ and $\nu>\frac{1}{q(H-1)+1}$, the result follows.

This proves that, by letting $m$ and then $n$ go to infinity, $A^{\left(
m,n\right) }$ converges in $L^{2}$ to $\int_{0}^{1}f\left( B_{s}\right)
dZ_{s}^{(q)}.$

\vskip0.5cm

\textit{Study of the term $B^{\left( m,n\right) }$: } We prove that
\begin{equation}
\lim_{n\rightarrow \infty }\sup_{m }E\left|
B^{(m,n)}\right| ^{2}=0.  \label{a2a}
\end{equation}%
We have, using the product formula (\ref{multiplication}) for multiple stochastic integrals,
\begin{eqnarray}
E\left| B^{(m,n)}\right| ^{2}
&=&2^{2m(q-1)}\sum_{j=1}^{2^{n}}\sum_{i=(j-1)2^{m-n}+1}^{j2^{m-n}}\sum_{j^{%
\prime }=1}^{2^{n}}\sum_{i^{\prime }=(j^{\prime }-1)2^{m-n}+1}^{j^{\prime
}2^{m-n}}\sum_{l=0}^{q}\frac{l!}{q!^{2}}\binom{q}{l}^{2}
\notag \\
&&\times b_{l}^{(m,n)}\langle \delta _{i2^{-m}},\delta _{i^{\prime
}2^{-m}}\rangle _{\EuFrak H}^{l},  \label{a4}
\end{eqnarray}%
where%
\begin{equation}
b_{l}^{(m,n)}=E\left( \Delta _{i,j}^{m,n}f(B)\Delta _{i^{\prime },j^{\prime
}}^{m,n}f(B)I_{2(q-l)}\left( \delta _{i2^{-m}}^{\otimes (q-l)}\widetilde{%
\otimes }\delta _{i^{\prime }2^{-m}}^{\otimes (q-l)}\right) \right) .
\label{bel}
\end{equation}
By (\ref{dual}) and (\ref{der-multiple}), we obtain that $b_l^{(m,n)}$ is equal to
\begin{eqnarray*}
&&E\left\langle D^{2(q-l)}\left( \Delta
_{i,j}^{m,n}f(B)\Delta _{i^{\prime },j^{\prime }}^{m,n}f(B)\right)
,\delta _{i2^{-m}}^{\otimes (q-l)}\widetilde{\otimes }%
\delta _{i^{\prime }2^{-m}}^{\otimes (q-l)} \right\rangle _{\EuFrak %
H^{\otimes 2(q-l)}} \\
&=&\sum_{a=0}^{2q-2l}\binom{2q-2l}{a}\left\langle E\left( \left(
f^{(a)}(B_{(i-1)2^{-m}})\varepsilon _{(i-1)2^{-m}}^{\otimes
a}-f^{(a)}(B_{(j-1)2^{-n}})\varepsilon _{(j-1)2^{-n}}^{\otimes a}\right)
\widetilde{\otimes} \right. \right.  \\
&&\left. \left.
\left( f^{(2q-2l-a)}(B_{(i^{\prime }-1)2^{-m}})\varepsilon
_{(i^{\prime }-1)2^{-m}}^{\otimes b}-f^{(2q-2l-a)}(B_{(j^{\prime
}-1)2^{-n}})\varepsilon _{(j^{\prime }-1)2^{-m}}^{\otimes b}\right)
\right) ,\delta _{i2^{-m}}^{\otimes (q-l)}\widetilde{\otimes} \delta
_{i^{\prime }2^{-m}}^{\otimes (q-l)}\right\rangle _{\EuFrak
H^{\otimes 2(q-l)}}.
\end{eqnarray*}
The term in (\ref{a4}) corresponding to $l=q$ can be estimated by%
\begin{equation*}
\frac{1}{q!}2^{2m(q-1)}\sup_{|x-y|\leq 2^{-n}}E\left|
f(B_{x})-f(B_{y})\right| ^{2}\beta _{q,m},
\end{equation*}%
where $\beta _{q,m}$ has been introduced in (\ref{beta}).
So it converges to
zero as $n$ tends to infinity, uniformly in $m$, because, by
(\ref{aux4new}) and using that $H>1-\frac{1}{2q}$, we have
\begin{equation*}
\sup_{m}2^{2m(q-1)}\beta _{q,m}<\infty .
\end{equation*}%
In order to handle the terms with $0\leq l\leq q-1$, we make the
decomposition
\begin{equation}
\left| b_{l}^{(m,n)}\right| \leq \sum_{a=0}^{2q-2l}\binom{2q-2l}{a}\sum_{h=1}^{4}B_{h},  \label{dec}
\end{equation}
where
\begin{eqnarray}
B_{1} &=&E\left| \Delta _{i,j}^{m,n}f(B)\Delta _{i^{\prime },j^{\prime
}}^{m,n}f(B)\right| \left\langle \varepsilon _{(i-1)2^{-m}}^{\otimes a}\ \widetilde{\otimes}
\varepsilon _{(i^{\prime }-1)2^{-m}}^{\otimes (2q-2l-a)},\delta _{i2^{-m}}^{\otimes
(q-l)}\widetilde{\otimes }\delta _{i^{\prime }2^{-m}}^{\otimes (q-l)}\right\rangle
_{\EuFrak H^{\otimes 2(q-l)}},  \notag \\
B_{2} &=&E\left| f^{(a)}(B_{(j-1)2^{-n}})\Delta _{i^{\prime },j^{\prime
}}^{m,n}f(B)\right|   \notag \\
&&\times \left\langle \left( \varepsilon _{(i-1)2^{-m}}^{\otimes a}-\varepsilon
_{(j-1)2^{-n}}^{\otimes a}\right)\widetilde{\otimes} \varepsilon _{(i^{\prime
}-1)2^{-m}}^{\otimes (2q-2l-a)},\delta _{i2^{-m}}^{\otimes (q-l)}\widetilde{\otimes }%
\delta _{i^{\prime }2^{-m}}^{\otimes (q-l)}\right\rangle _{\EuFrak H^{\otimes
2(q-l)}},  \notag \\
B_{3} &=&E\left| \Delta _{i,j}^{m,n}f(B)f^{(2q-2l-a)}(B_{(j^{\prime
}-1)2^{-n}})\right|   \notag \\
&&\times \left\langle \varepsilon _{(i-1)2^{-m}}^{\otimes a}\widetilde{\otimes} \left(
\varepsilon _{(i^{\prime }-1)2^{-m}}^{\otimes (2q-2l-a)}-\varepsilon _{(j^{\prime
}-1)2^{-n}}^{\otimes (2q-2l-a)}\right) ,\delta _{i2^{-m}}^{\otimes (q-l)}\widetilde{%
\otimes }\delta _{i^{\prime }2^{-m}}^{\otimes (q-l)}\right\rangle _{\EuFrak %
H^{\otimes 2(q-l)}},  \notag \\
B_{4} &=&E\left| f^{(a)}(B_{(j-1)2^{-n}})f^{(2q-2l-a)}(B_{(j^{\prime
}-1)2^{-n}})\right|   \notag \\
&&\times \left\langle \left( \varepsilon _{(i-1)2^{-m}}^{\otimes a}-\varepsilon
_{(j-1)2^{-n}}^{\otimes a}\right) \widetilde{\otimes} \left( \varepsilon _{(i^{\prime
}-1)2^{-m}}^{\otimes (2q-2l-a)}-\varepsilon _{(j^{\prime }-1)2^{-n}}^{\otimes
(2q-2l-a)}\right) ,\delta _{i2^{-m}}^{\otimes (q-l)}\widetilde{\otimes }\delta
_{i^{\prime }2^{-m}}^{\otimes (q-l)}\right\rangle _{\EuFrak H^{\otimes 2(q-l)}}. \notag \\
\label{bi}
\end{eqnarray}
By using (\ref{aux5}) and the conditions imposed on the
function $f$, one can bound the terms $B_{1}$, $B_{2}$ and $B_{3}$
as follows:
\begin{equation*}
|B_{1}|\leq c(q,f,H)\sup_{|x-y|\leq \frac{1}{2^{n}},0\leq a\leq
2q}E\left| f^{(a)}(B_{x})-f^{(a)}(B_{y})\right| ^{2}2^{-2m(q-l)},
\end{equation*}

\begin{equation*}
|B_{2}|+|B_{3}|\leq c(q,f,H)\sup_{|x-y|\leq \frac{1}{2^{n}},0\leq
a\leq 2q}E\left| f^{(2q-2l-a)}(B_{x})-f^{(2q-2l-a)}(B_{y})\right| 2^{-2m(q-l)},
\end{equation*}
and, by using (\ref{boundHgrand}), we obtain that
\begin{equation*}
|B_{4}|\leq c(q,f,H)2^{-n\frac{q-1}{q}-2m(q-l)}.
\end{equation*}
By setting
$$
R_n=\frac{1}{q!}\sup_{|x-y|\leq 2^{-n}}E\left|
f(B_{x})-f(B_{y})\right| ^{2}\sup_{m}2^{2m(q-1)}\beta _{q,m},
$$
we can finally write, by the estimate (\ref{aux4new}),
\begin{eqnarray*}
&&E\left| B^{(m,n)}\right| ^{2}\\ &\leq &R_{n}+c(H,f,q)2^{2m(q-1)}\left(
\sup_{|x-y|\leq \frac{1}{2^{n}},0\leq a\leq 2q}\left|
f^{(2q-2l-a)}(B_{x})-f^{(2q-2l-a)}(B_{y})\right| + (2^{-n})^{\frac{q-1}{q}} \right)  \\
&&\times \sum_{j=1}^{2^{n}}\sum_{i=(j-1)2^{m-n}+1}^{j2^{m-n}}\sum_{j^{\prime
}=1}^{2^{n}}\sum_{i^{\prime }=(j^{\prime }-1)2^{m-n}+1}^{j^{\prime
}2^{m-n}}\sum_{l=0}^{q-1}2^{-2m(q-l)}\langle \delta
_{i2^{-m}},\delta _{i^{\prime }2^{-m}}\rangle _{\EuFrak H}^{l} \\
&\leq &R_{n}+c(H,f,q)2^{2m(q-1)}\left( \sup_{|x-y|\leq \frac{1}{2^{n}},0\leq
a\leq 2q}\left| f^{(2q-2l-a)}(B_{x})-f^{(2q-2l-a)}(B_{y})\right| +(2^{-n})^{\frac{q-1}{q}%
}\right)  \\
&&\times \sum_{l=0}^{q-1}2^{-2m(q-l)}\sum_{i,j=0}^{2^{m}}\langle
\delta _{i2^{-m}},\delta _{i^{\prime }2^{-m}}\rangle _{\EuFrak H}^{l} \\
&\leq &R_{n}+c(H,f,q)\left( \sup_{|x-y|\leq \frac{1}{2^{n}},0\leq a\leq
2q}\left| f^{(2q-2l-a)}(B_{x})-f^{(2q-2l-a)}(B_{y})\right| +(2^{-n})^{\frac{q-1}{q}%
}\right)
\end{eqnarray*}%
and this converges to zero due to the continuity of $B$ and since
$q>1$. \qed
\vskip0.3cm

\subsection{Proof of Theorem \ref{main} in the case $\frac{1}{2q}<H\leq 1-%
\frac{1}{2q}$: the weighted central limit theorem}\label{subs}

Suppose first that $\frac{1}{2q}<H<1-\frac{1}{2q}$.
%%% debut
We study the convergence in law of the sequence
  $V_{2,n}^{(q)}(f)=2^{-\frac{n}{2}}\,V_{n}^{(q)}(f)$.
We fix two integers $m\geq n$, and decompose this sequence as follows:
\begin{equation*}
V_{2,m}^{(q)}(f)=A^{(m,n)}+B^{(m,n)},
\end{equation*}%
where
\begin{equation*}
A^{(m,n)}=2^{-\frac{m}{2}}\sum_{j=1}^{2^{n}}f\left( B_{(j-1)2^{-n}}\right)
\sum_{i=(j-1)2^{m-n}+1}^{j2^{m-n}}H_{q}\left( 2^{mH}
%%% debut
\Delta B_{i2^{-m}}
%%% fin
\right) ,
\end{equation*}%
and
\begin{equation*}
B^{\left( m,n\right) }=\frac{1}{q!}\,2^{m(Hq-\frac{1}{2})}\sum_{j=1}^{2^{n}}%
\sum_{i=(j-1)2^{m-n}+1}^{j2^{m-n}}\Delta _{i,j}^{m,n}f(B)I_{q}\left( \delta
_{i2^{-m}}^{\otimes q}\right) ,
\end{equation*}%
and where as before we make use of the notation $\Delta
_{i,j}^{m,n}f(B)=f(B_{\left( i-1\right) 2^{-m}})-f(B_{\left( j-1\right)
2^{-n}})$. \

Let us first consider the term $A^{(m,n)}$. From
Theorem 1 in Breuer
and Major \cite{BM}, and taking into account that  $H<1-\frac{1}{2q}$,
it follows that the random vector
$$
\left(B,
2^{-\frac{m}2}\sum_{i=(j-1)2^{m-n}+1}^{j2^{m-n}}H_q(2^{mH}\Delta B_{i2^{-m}});\quad j=1,\ldots,2^n
\right)
$$
converges in law, as $m\to\infty$, to
$$
\big(B,\sigma_{H,q}\Delta W_{j2^{-n}};\quad j=1,\ldots,2^n\big)
$$
where
$\sigma_{H,q}$
is the constant defined by (\ref{hsmall}) and
$W$  is a standard Brownian motion independent of $B$
(the independence is a consequence of the central limit theorem for multiple stochastic integrals proved
in Peccati and Tudor \cite{PT}).
Since
\begin{equation*}
\sum_{j=1}^{2^{n}}f\left( B_{(j-1)2^{-n}}\right) \,\Delta W_{j2^{-n}}
\end{equation*}%
converges in $L^{2}$ as $n\rightarrow \infty $ to the It\^{o} integral $%
\int_{0}^{1}f(B_{s})dW_{s}$ we conclude that, by letting $m\rightarrow
\infty $ and then $n\rightarrow \infty $, we have
\begin{equation*}
\big(B,A^{(m,n)}\big)\,\,\,{\overset{\mathrm{Law}}{\longrightarrow }}\,\,
\left(B,\sigma_{H,q}%
\int_{0}^{1}f(B_{s})dW_{s}\right).
\end{equation*}%
Then it suffices to show that
\begin{equation}
\lim_{n\rightarrow \infty }  \sup_{m\rightarrow \infty }E\left|
B^{(m,n)}\right| ^{2}=0.  \label{a1}
\end{equation}%
We have, as in (\ref{a4}),
\begin{eqnarray}
E\left| B^{(m,n)}\right| ^{2}
&=&2^{m(2Hq-1)}\sum_{j=1}^{2^{n}}\sum_{i=(j-1)2^{m-n}+1}^{j2^{m-n}}\sum_{j^{%
\prime }=1}^{2^{n}}\sum_{i^{\prime }=(j^{\prime }-1)2^{m-n}+1}^{j^{\prime
}2^{m-n}}\sum_{l=0}^{q}\frac{l!}{q!^{2}}\binom{q}{l}^{2}
\notag \\
&&\times b_{l}^{(m,n)}\langle \delta _{i2^{-m}},\delta _{i^{\prime
}2^{-m}}\rangle _{\EuFrak H}^{l},  \label{a2}
\end{eqnarray}%
where $b_{l}^{(m,n)}$ has been defined in (\ref{bel}). The term in (\ref{a2}%
) corresponding to $l=q$ can be estimated by%
\begin{equation*}
\frac{1}{q!}2^{m(2Hq-1)}\sup_{|x-y|\leq 2^{-n}}E\left|
f(B_{x})-f(B_{y})\right| ^{2}\beta _{q,m},
\end{equation*}%
which converges to zero as $n$ tends to infinity, uniformly in $m$, because
by (\ref{aux4}) and using that $H<1-\frac{1}{2q}$, we have
\begin{equation*}
\sup_{m}2^{m(2Hq-1)}\beta _{q,m}<\infty .
\end{equation*}
In order to handle the terms with $0\leq l\leq q-1$, we will distinguish two
different cases, depending on the value of \ $H$.

\vskip0.3cm

\textit{Case }$H<1/2$. \ Suppose $0\leq l\leq q-1$. By (\ref{aux2}), we can majorize $b_{l}^{(m,n)}$ as follows:
\begin{equation*}
|b_{l}^{(m,n)}|\leq C2^{-4Hm(q-l)}.
\end{equation*}
As a consequence, applying again (\ref{aux4}), the
corresponding term in (\ref{a2}) is bounded by
\begin{equation*}
C2^{m(2Hq-1)}2^{-4Hm(q-l)}\beta _{l,m}\leq C2^{2mH(l-q)},
\end{equation*}%
which converges to zero as $m$ tends to infinity because $l<q$.

\vskip0.3cm

\textit{Case }$H>1/2$. \ Suppose $0\leq l\leq q-1$. By (\ref{aux5}), we get the estimate
\begin{equation*}
|b_{l}^{(m,n)}|\leq C2^{-2m(q-l)}.
\end{equation*}
As a consequence, applying again (\ref{aux4}), the
corresponding term in (\ref{a2}) is bounded by
\begin{equation*}
C2^{m(2Hq-1)}2^{-2m(q-l)}\beta _{l,m}.
\end{equation*}%
If \ $H<1-\frac{1}{2l}$, applying (\ref{aux4}), this is
bounded by $C2^{m(2H(q-l)-2(q-l))}$, which converges to zero as $m$ tends to
infinity because $H<1$ and $l<q$. In the case $H=1-\frac{1}{2l}$, applying
(\ref{aux6}), we get the estimate $Cm2^{m(2H(q-l)-2(q-l))}
$, which converges to zero as $m$ tends to infinity because $H<1$ and $%
l<q$. In the case $H>1-\frac{1}{2l}$, we apply (\ref{aux5})
and we get the estimate $C2^{m(2H2+1-2q)}$, which converges to zero as $m$
tends to infinity because $H<1-\frac{1}{2q}$.

\bigskip

The proof in the case $H=1-\frac{1}{2q}$ is similar. The convergence of the
term \ $A^{(m,n)}$ is obtained by applying Theorem 1' in Breuer and Major
(1983), and the convergence to zero in $L^{2}$ of the term $B^{(m,n)}$
follows the same lines as before.

%%% fin

\vskip0.3cm
%%% debut
\subsection{Proof of Proposition \ref{joint}}
We proceed as in Section \ref{subs}.
For $p=2,\dots,q$, we set
  $V_{2,n}^{(p)}(f)=2^{-\frac{n}{2}}\,V_{n}^{(p)}(f)$.
We fix two integers $m\geq n$, and decompose this sequence as follows:
\begin{equation*}
V_{2,m}^{(p)}(f)=A_p^{(m,n)}+B_p^{(m,n)},
\end{equation*}%
where
\begin{equation*}
A_p^{(m,n)}=2^{-\frac{m}{2}}\sum_{j=1}^{2^{n}}f\left( B_{(j-1)2^{-n}}\right)
\sum_{i=(j-1)2^{m-n}+1}^{j2^{m-n}}H_{p}\left( 2^{mH}\Delta B_{i2^{-m}}\right) ,
\end{equation*}%
and
\begin{equation*}
B_p^{\left( m,n\right) }=\frac{1}{p!}\,2^{m(Hp-\frac{1}{2})}\sum_{j=1}^{2^{n}}%
\sum_{i=(j-1)2^{m-n}+1}^{j2^{m-n}}\Delta _{i,j}^{m,n}f(B)I_{p}\left( \delta
_{i2^{-m}}^{\otimes p}\right) ,
\end{equation*}%
and where as before we make use of the notation $\Delta
_{i,j}^{m,n}f(B)=f(B_{\left( i-1\right) 2^{-m}})-f(B_{\left( j-1\right)
2^{-n}})$.

Let us first consider the term $A_p^{(m,n)}$. We claim that   the random vector
\begin{equation*}
\left( B,\left\{2^{-\frac{m}{2}}\sum_{i=(j-1)2^{m-n}+1}^{j2^{m-n}}H_{p}\left(
2^{mH}\Delta B_{i2^{-m}}\right) ;\quad j=1,\ldots ,2^{n}\right\}_{2\leq p\leq q}\right)
\end{equation*}%
converges in law, as $m\rightarrow \infty $, to
\begin{equation*}
\left( B,\{\sigma_{H,p}\Delta W^{(p)}_{j2^{-n}};\quad j=1,\ldots ,2^{n}\}_{2\leq p\leq q}\right)
\end{equation*}%
where $(W^{(2)}, \dots, W^{(q)})$ is a $(q-1)$-dimensional standard Brownian motion independent of $B$
%%% debut
and the $\sigma_{H,p}$'s are given by (\ref{hsmall}).
%%% fin
Indeed,
the convergence in law of each component follows from  Theorem 1 in Breuer
and Major \cite{BM}, taking into account that  $H<\frac34\leq 1-\frac{1}{2q}$.  The joint convergence
and the fact that the processes
$W^{(p)}$ for $p=2,\dots,q$ are independent (and also independent of $B$)
%%% debut
is a direct application of
%%% fin
the central limit theorem
for multiple stochastic integrals proved
in Peccati and Tudor \cite{PT}.
%%% debut
%Indeed, the orthogonality of Herite polynomials of different order with
%respect to the Gaussian distribution implies
%the independence of the components of the limit vector.
%
%We refer to Theorem 1 in Breuer
%and Major \cite{BM} for the explicit value of the constants  $\sigma_{H,p}$.
%%% fin

Since, for any $p=2,\ldots,q$, the quantity
\begin{equation*}
\sum_{j=1}^{2^{n}}f\left( B_{(j-1)2^{-n}}\right) \,\Delta W^{(p)}_{j2^{-n}}
\end{equation*}%
converges in $L^{2}$ as $n\rightarrow \infty $ to the It\^{o} integral $%
\int_{0}^{1}f(B_{s})dW^{(p)}_{s}$, we conclude that, by letting $m\rightarrow
\infty $ and then $n\rightarrow \infty $, we have
\begin{equation*}
\left(
%%% debut
B,
%%% fin
A_2^{(m,n)}, \dots, A_q^{(m,n)}  \right)\,\,\,{\overset{\mathrm{Law}}{\longrightarrow }}\,\,
\left(
%%% debut
B,
%%% fin
\sigma_{H,2}
\int_{0}^{1}f(B_{s})dW^{(2)}_{s}, \dots, \sigma_{H,q}%
\int_{0}^{1}f(B_{s})dW^{(q)}_{s}  \right).
\end{equation*}%
On the other hand, and because $H\in(\frac14,\frac34)$ (implying that
$H\in(\frac1{2p},1-\frac{1}{2p})$), we have shown in Section \ref{subs} that
$$
\lim_{n\rightarrow \infty }  \sup_{m\rightarrow \infty }E\left|
B_p^{(m,n)}\right| ^{2}=0
$$
for all $p=2,\ldots,q$.
This finishes the proof of Proposition \ref{joint}.

%%% fin
\subsection{Proof of Corollary \ref{complete-picture}}

For any integer $q\geq 2$, we have
\begin{eqnarray*}
\left( 2^{nH}\Delta B_{k2^{-n}}\right) ^{q}-\mu _{q} =\sum_{p=1}^{q}
\binom{q}{p}\mu _{q-p}2^{Hnp}I_{p}(\delta _{k2^{-n}}^{\otimes p})
=\sum_{p=1}^{q}p!\binom{q}{p}\mu _{q-p}H_{p}\left( 2^{nH}\Delta
B_{k2^{-n}}\right).
\end{eqnarray*}
Indeed,
the $p$th kernel in the chaos representation of \ $\left(
2^{nH}\Delta B_{k2^{-n}}\right) ^{q}$ is
\begin{equation*}
\frac{1}{p!}E(D^{p}\left( 2^{nH}\Delta B_{k2^{-n}}\right) ^{q})=
\binom{q}{p}2^{nHp}\mu _{q-p}\delta _{k2^{-n}}^{\otimes p}.
\end{equation*}

Suppose first that $q$ is odd and $H>\frac{1}{2}$. In this case, we have
\begin{equation*}
2^{-nH}\sum_{k=1}^{2^{n}}f(B_{(k-1)2^{-n}})(2^{nH}\Delta
B_{k2^{-n}})^{q}=\sum_{p=1}^{q}p!\binom{q}{p}\mu _{q-p}2^{-nH}V_{n}^{(p)}(f).
\end{equation*}%
The term with $p=1$ converges in $L^{2}$ to $q\mu
_{q-1}\int_{0}^{1}f(B_{s})dB_{s}$. For $p\geq 2$, the limit in
$L^{2}$ is zero. Indeed, if $H\le 1-\frac{1}{2p}$, then $E\left(
V_{n}^{(p)}(f)^{2}\right) $ is bounded by a constant times $n2^{n}$
by Proposition \ref{p1}. If $H>1-\frac{1}{2p}$, then $E\left(
V_{n}^{(p)}(f)^{2}\right)$ is bounded by a constant times
$2^{-n2(1-H)p+2n}$ by (\ref{t3}),
with $-2(1-H)p+2-2H=(1-H)(2-2p)<0$.

Suppose now that $q$ is even. Then
\begin{equation*}
2^{2nH-n}\sum_{k=1}^{2^{n}}f(B_{(k-1)2^{-n}})\left[ (2^{nH}\Delta
B_{k2^{-n}})^{q}-\mu _{q}\right] =2^{2nH-n}\sum_{p=2}^{q}p!\binom{q}{p}\mu
_{q-p}V_{n}^{(p)}(f)\ .
\end{equation*}

If \ $H<\frac{1}{4}$, by (\ref{t1}), one has that $2^{2nH-n}\times 2\binom{q}{%
2}\mu _{q-2}V_{n}^{(2)}(f)$ converges in $L^{2}$, as $n\to\infty$, to
$\frac{1}{4}\binom{q}{2}\mu _{q-2}\int_{0}^{1}f^{\prime \prime
}(B_{s})ds$. On the other hand, for $p\geq 4$, $2^{2nH-n}V_{n}^{(p)}(f)$
converges to zero in $L^{2}$. Indeed, if $H<\frac{1}{2p}$, then
$E\left( V_{n}^{(p)}(f)^{2}\right) =O(2^{n(-2Hp+2)})$ by (\ref%
{r1}) with $%
-2Hp+2+4H-2<0$. If $H\geq \frac{1}{2p}$, then $E\left( V_{n}^{(p)}(f)^{2}\right) =O(2^{n})$
by \ (\ref{r1a}) with $4H-1<0$. Therefore (\ref%
{c2}) holds.

In the case $\frac{1}{4}<H<\frac{3}{4}$,
Proposition \ref{joint} implies that  the vector
\[
\left(B,2^{-n/2}V_{n}^{(2)}(f), \dots, 2^{-n/2}V_{n}^{(q)}(f) \right)
\]
 converges in law
to
\[
\left( B,\sigma _{H,2}\int_{0}^{1}f(B_{s})dW_{s}^{(2)},   \ldots,
\sigma _{H,q}\int_{0}^{1}f(B_{s})dW_{s}^{(q)}  \right),
\]
where $(W^{(2)},\ldots,W^{(q)})$ is a $(q-1)$-dimensional standard Brownian motion independent of $B$
and the $\sigma_{H,p}$'s, $2\leq p\leq q$, are given by (\ref{hsmall}).
This implies the convergence \ (\ref{c3}). The proof of
(\ref{c3a}) is analogous
%%% debut
(with an adequate version of Proposition \ref{joint}).
%%% fin

The convergence (\ref{c3a1/4}) is obtained by similar arguments using the limit result
(\ref{conjecture}) in the critical case $H=\frac 14$, $p=2$.

Finally, consider the case \ $H>\frac{3}{4}$. For $p=2$,
$2^{n-2Hn}V_{n}^{(2)}(f)$ converges in $L^{2}$ to
$
%%% debut
%2
%%% fin
\int_{0}^{1}f(B_{s})dZ_{s}^{(2)}$ by (\ref{t3}). If $%
p\geq 4$, then $2^{n-2Hn}V_{n}^{(p)}(f)$ $\ $converges in $L^{2}$ to zero
because, again by (\ref{t3}), one has $E\left(
V_{n}^{(p)}(f)^{2}\right) =O(2^{n(2-2(1-H)p)})$.

\subsection{Proof of Theorem \ref{ito-thm}}

We   can assume $H<\frac{1}{2}$, the case where $H\geq \frac{1}{2}$
being straightforward. By a Taylor's formula, we have
\begin{eqnarray}
f(B_{1}) &=&f(0)+\frac{1}{2}\sum_{k=1}^{2^{n}}\big(f^{\prime
}(B_{k2^{-n}})+f^{\prime }(B_{(k-1)2^{-n}})\big)\,\Delta B_{k2^{-n}}  \notag
\\
&&-\frac{1}{12}\sum_{k=1}^{2^{n}}f^{(3)}(B_{(k-1)2^{-n}})\big(\Delta
B_{k2^{-n}}\big)^{3}-\frac{1}{24}\sum_{k=1}^{2^{n}}f^{(4)}(B_{(k-1)2^{-n}})%
\big(\Delta B_{k2^{-n}}\big)^{4}  \notag \\
&&-\frac{1}{80}\sum_{k=1}^{2^{n}}f^{(5)}(B_{(k-1)2^{-n}})\big(\Delta
B_{k2^{-n}}\big)^{5}+R_{n},  \label{taylor}
\end{eqnarray}%
with $R_{n}$ converging towards $0$ in probability as $n\rightarrow \infty $%
, because $H>1/6$.  We can expand the monomials $x^{m}$, $m=2,3,4,5$, in
terms of the Hermite polynomials:
\begin{eqnarray*}
x^{2} &=&2\,H_{2}(x)+1 \\
x^{3} &=&6\,H_{3}(x)+3\,H_{1}(x) \\
x^{4} &=&24\,H_{4}(x)+12\,H_{2}(x)+3 \\
x^{5} &=&120\,H_{5}(x)+60\,H_{3}(x)+15\,H_{1}(x).
\end{eqnarray*}%
In this way we obtain
\begin{eqnarray}
&&\sum_{k=1}^{2^{n}}f^{(3)}(B_{(k-1)2^{-n}})\left( \Delta B_{k2^{-n}}\right)
^{3} =6\times 2^{-3Hn}V_{n}^{(3)}(f^{(3)})+3\times 2^{-2Hn}V_{n}^{(1)}(f^{(3)}),
\label{h1} \\
&& \sum_{k=1}^{2^{n}}f^{(4)}(B_{(k-1)2^{-n}})\left( \Delta B_{k2^{-n}}\right)
^{4} =24\times 2^{-4Hn}V_{n}^{(4)}(f^{(4)})  \notag  \\
 &&  \qquad \qquad  +12\times
2^{-4Hn}V_{n}^{(2)}(f^{(4)})+3\times
2^{-4Hn}\sum_{k=1}^{2^{n}}f^{(4)}(B_{(k-1)2^{-n}}),  \label{h2} \\
&& \sum_{k=1}^{2^{n}}f^{(5)}(B_{(k-1)2^{-n}})\left( \Delta B_{k2^{-n}}\right)
^{5} =120\times 2^{-5Hn}V_{n}^{(5)}(f^{(5)})\notag  \\
&&  \qquad \qquad +60\times
2^{-5Hn}V_{n}^{(3)}(f^{(5)})+15\times 2^{-4Hn}V_{n}^{(1)}(f^{(5)}).
\label{h3}
\end{eqnarray}%
By (\ref{r1a}) and using that  $H>\frac{1}{6}$, we have $E\left(
V_{n}^{(3)}(f^{(3)})^{2}\right) \leq C2^{n}$ and $E\left(
V_{n}^{(3)}(f^{(5)})^{2}\right) \leq C2^{n}$. As a consequence, the
first summand in\ (\ref{h1}) and the second one in (\ref{h3}) converge  to
zero in $L^{2}$ as $n$ tends to infinity. \ Also, by (\ref{r1a}), $E\left(
V_{n}^{(4)}(f^{(4)})^{2}\right) \leq C2^{n}$ and $E\left(
V_{n}^{(5)}(f^{(5)})^{2}\right) \leq C2^{n}$. Hence, the first summand in (\ref%
{h2}) and the first summand in (\ref{h3}) converge  to zero in $L^{2}$ as $n$
tends to infinity.   If $\frac{1}{6}<H<\frac{1}{4}$, (\ref{r1}) implies  $%
E\left( V_{n}^{(2)}(f^{(4)})^{2}\right) \leq C2^{n(-4H+2)})$, so that $%
2^{-4Hn}V_{n}^{(2)}(f^{(4)})$ converges to zero in $L^{2}$ as $n$
tends to infinity. If $\frac{1}{4}\leq H<\frac{1}{2}$, (\ref{r1a}) implies
$E\left( V_{n}^{(2)}(f)^{2}\right) \leq C2^{n}$ so that
$2^{-4Hn}V_{n}^{(2)}(f^{(4)})$ converges to zero in $L^{2}$ as $n$
tends to infinity.

Moreover, using the following identity, valid for regular functions $h:\mathbb{R}%
\rightarrow \mathbb{R}$:
\begin{equation*}
\sum_{k=1}^{2^{n}}h^{\prime }(B_{(k-1)2^{-n}})\,\Delta
B_{k2^{-n}}=h(B_{1})-h(0)-\frac{1}{2}\sum_{k=1}^{2^{n}}h^{\prime \prime
}(B_{\theta _{k2^{-n}}})\left( \Delta B_{k2^{-n}}\right) ^{2}
\end{equation*}%
for some $\theta _{k2^{-n}}$ lying between $(k-1)2^{-n}$ and $k2^{-n}$, we
deduce that $2^{-4Hn}V_{n}^{(1)}(f^{(5)})$ tends to zero, because $H>\frac{1%
}{6}$. In the same way,  we have
\begin{eqnarray*}
2^{-2Hn}V_{n}^{(1)}(f^{(3)}) &=&-\frac{1}{2}\,2^{-2Hn}%
\sum_{k=1}^{2^{n}}f^{(4)}(B_{(k-1)2^{-n}})\left( \Delta B_{k2^{-n}}\right)
^{2} \\
&&-\frac{1}{6}\,2^{-2Hn}\sum_{k=1}^{2^{n}}f^{(5)}(B_{(k-1)2^{-n}})\left(
\Delta B_{k2^{-n}}\right) ^{3}+o(1).
\end{eqnarray*}%
We have obtained
\begin{eqnarray*}
f(B_{1}) &=&f(0)+\frac{1}{2}\sum_{k=1}^{2^{n}}\left( f^{\prime
}(B_{k2^{-n}})+f^{\prime }(B_{(k-1)2^{-n}})\,\right) \Delta B_{k2^{-n}} \\
&&+\frac{1}{4}\times 2^{-4Hn}\sum_{k=1}^{2^{n}}f^{(4)}(B_{(k-1)2^{-n}})H_{2}%
\left( 2^{nH}\Delta B_{k2^{-n}}\right) \  \\
&&-\frac{1}{24}\times 2^{-2Hn}\sum_{k=1}^{2^{n}}f^{(5)}(B_{(k-1)2^{-n}})\left(
\Delta B_{k2^{-n}}\right) ^{3}+o(1).
\end{eqnarray*}%
As before $\,2^{-4Hn}V_{n}^{(2)}(f^{(4)})$ converges to zero in $L^{2}$.
Finally, by (\ref{pomfbim}),
$$
2^{-2Hn}\sum_{k=1}^{2^{n}}f^{(5)}(B_{(k-1)2^{-n}})\left( \Delta
B_{k2^{-n}}\right) ^{3}
$$
also converges to zero. This
completes the proof. \qed

\vskip0.5cm
\noindent
{\bf Acknowledgments}.
%%% debut
We are grateful to Jean-Christophe Breton and Nabil
Kazi-Tani for helpful remarks. We also wish to thank the anonymous referee for
his/her very careful reading.
%%% fin

\vskip 0.5cm

%%% j'ai mis à jour la biblio

\end{document}